\documentclass[10pt,reqno]{amsart}
\usepackage{hyperref}
\usepackage{float}
\usepackage{amscd,amssymb,amsmath,latexsym,bm}
\usepackage{mathtools}
\usepackage[mathcal,mathscr]{euscript}
\usepackage{lipsum}
\usepackage{amsfonts}
\usepackage{graphicx}
\usepackage{epstopdf}
\usepackage{algorithmic}
\usepackage{hyperref}
\usepackage{xcolor}
\usepackage{upgreek}
\usepackage{enumerate}
\usepackage{ulem}
\usepackage{stmaryrd}
\usepackage{gensymb}			
\usepackage[utf8]{inputenc}
\usepackage{mathtools}
\newcommand*{\twoheadrightarrowtail}{\mathrel{\rightarrowtail\kern-1.9ex\twoheadrightarrow}}

\usepackage{tikz-cd}
\usepackage{quiver}
 \usetikzlibrary{shapes,arrows}
 
\newtheorem{theorem}{Theorem}[section]

\newtheorem{corollary}[theorem]{Corollary}
\newtheorem{proposition}[theorem]{Proposition}
\newtheorem{remark}[theorem]{Remark}
\newtheorem{definition}[theorem]{Definition}
\newtheorem{example}[theorem]{Example}

\numberwithin{equation}{section}
\numberwithin{figure}{section}

\newcommand{\AM}{{\mathbb A}}

\newcommand{\CM}{{\mathbb C}}
\newcommand{\NM}{{\mathbb N}}

\newcommand{\RM}{{\mathbb R}}

\newcommand{\ZM}{{\mathbb Z}}
\newcommand{\PM}{{\mathbb P}}

\newcommand{\EM}{{\mathbb E}}

\newcommand{\Aa}{{\mathcal A}}

\newcommand{\Ii}{{\mathcal I}}
\newcommand{\Pp}{{\mathcal P}}

\newcommand{\Bb}{{\mathcal B}}
\newcommand{\Dd}{{\mathcal D}}

\newcommand{\Ww}{{\mathcal W}}

\newcommand{\Vv}{{\mathcal V}}
\newcommand{\Ss}{{\mathcal S}}
\newcommand{\Oo}{{\mathcal O}}
\newcommand{\Tt}{{\mathcal T}}

\newcommand{\Mm}{{\mathcal M}}
\newcommand{\Cc}{{\mathcal C}}

\newcommand{\Jj}{{\mathcal J}}

\newcommand{\Qq}{{\mathcal Q}}
\newcommand{\Kk}{{\mathcal K}}

\begin{document}

\title[Operator Product States on Tensor Powers of $C^\ast$-Algebras]{Operator Product States on Tensor Powers of $C^\ast$-Algebras}

\author{Emil Prodan}

\address{Department of Physics and
\\ Department of Mathematical Sciences 
\\Yeshiva University 
\\New York, NY 10016, USA \\
\href{mailto:prodan@yu.edu}{prodan@yu.edu}}

\date{\today}

\begin{abstract} 
The program of matrix product states on tensor powers $\Aa^{\otimes \ZM}$ of $C^\ast$-algebras, initiated in Comm. Math. Phys. {\bf 144}, 443-490 (1992), is re-assessed in a context where $\Aa$ is a generic nuclear $C^\ast$-algebra. For any shift invariant state $\omega$, we demonstrate the existence of an order kernel ideal $\Kk_\omega$, whose quotient action reduces and factorizes the initial data $(\Aa^{\otimes \ZM}, \omega)$ to the tuple $(\Aa,\Bb_\omega = \Aa^{\otimes \NM^\times}/\Kk_\omega, \EM_\omega : \Aa \otimes \Bb_\omega \to \Bb_\omega, \bar \omega : \Bb_\omega \to \CM)$, where $\Bb_\omega$ is an operator system and  $\EM_\omega$ and $\bar \omega$ are unital and completely positive maps. Reciprocally, given a (input) tuple $(\Aa,\Ss,\EM,\phi)$ that shares similar attributes, we supply an algorithm that produces a shift-invariant state on $\Aa^{\otimes \ZM}$. We give sufficient conditions in which the so constructed states are ergodic and they reduce back to their input data. As examples, we formulate the input data that produces AKLT-type states, this time in the context of infinite site algebras, such as the group algebra of discrete amenable groups.
\end{abstract}

\thanks{This work was supported by the U.S. National Science Foundation through the grants DMR-1823800 and CMMI-2131760 and by U.S. Army Research Office through contract W911NF-23-1-0127.}

\maketitle


\setcounter{tocdepth}{3}

\section{Introduction and Main Statements}

Fannes,  Nachtergaele and Werner introduced and proved the following statement in their influential work \cite{FannesCMP1992}:

\begin{proposition}[\cite{FannesCMP1992}] Let $\Aa$ be a $C^\ast$-algebra with unit, and let $\omega$ be a shift-invariant state on $\Aa^{\otimes \ZM}$ {\rm (}$\simeq \Aa_L \otimes \Aa_R$, $\Aa_R = \Aa^{\otimes \NM^\times}${\rm )}.\footnote{In \cite{FannesCMP1992}, $\otimes$ refers exclusively to the minimal tensor product.} Then the following are equivalent:

\noindent (1) The set of functionals $\{\omega_x : \Aa_R \to \CM, \ \omega_x(a_R) = \omega(x \otimes a_R), \ x \in \Aa_L\}$, generate a finite dimensional linear subspace of $\Aa_R^\ast$.

\noindent (2) There are a finite dimensional vector space $\Bb$, a linear map $\EM$ from $\Aa$ to the space of linear maps over $\Bb$, an element $e \in \Bb$ and a linear functional $\rho \in \Bb^\ast$, such that $\rho\circ \EM_1 = \rho$, $\EM_1(e)=e$, and 
\begin{equation}\label{Eq:Facto0}
\omega(a_1 \otimes \cdots \otimes a_n) = \rho(e)^{-1} \rho \circ \EM_{a_1} \circ \cdots \circ \EM_{a_n}(e).
\end{equation}
\end{proposition} 

The states displaying this property with $\Bb$ a finite dimensional $C^\ast$-algebra were called $C^\ast$-finitely correlated states in \cite{FannesCMP1992}. They were shown to form a $\ast$-weakly dense convex subset of the set of translation invariant states on $\Aa^{\otimes \ZM}$, if $\Aa$ is a finite $C^\ast$-algebra.\footnote{In Proposition~2.6 of \cite{FannesCMP1992} that deals with this aspect, the space $\tilde \Bb$ defined there is finite dimensional only if $\Aa$ is finite dimensional. This detail that was omitted in \cite{FannesCMP1992}.}  In such cases, the identity~\eqref{Eq:Facto0} shows that these states factorize through the map $\EM$ and that the evaluation of such state on monomials is determined by a product of matrices. For this reason, these states are referred to as matrix-product states in the physics literature \cite{KlumperJPA1991,KlumperZPB1992,PeresQIC2007}. The work \cite{FannesCMP1992} sent the powerful message that any shift-invariant state on $\Aa^{\otimes \ZM}$ has arbitrarily close $\ast$-weak approximations that can be generated from the extremely simple data $(\Aa,\Bb,\EM,\rho,e)$. Furthermore, one can identify the ergodic states \cite{BratteliBook1}[Sec.~4.3] among the shift-invariant states by a simple examination of the spectral properties of $\EM$. These findings had a profound impact on the research on quantum spin chain systems.

In \cite{FannesCMP1992}, the authors pointed out that a fixed finitely correlated state $\omega$ has a minimal space $\Bb_\omega$ among all possible $\Bb$'s, which is uniquely determined by the state. One of our observations is that this space can be defined for any state on $\Aa^{\otimes \ZM}$ with $\Aa$ a nuclear $C^\ast$-algebra. Indeed, if we introduce what we call the entanglement kernel of a state 
\begin{equation}
\Kk_\omega = \bigcap_{x \in \Aa_L} {\rm Ker} \, \omega_x,
\end{equation}
then it is straightforward to see that $\Bb_\omega = \Aa_R/\Kk_\omega$ is exactly the minimal space mentioned in \cite{FannesCMP1992}, if the state happens to be finitely correlated. Furthermore, our second observation is that the map $\EM : \Aa \otimes \Bb_\omega \to \Bb_\omega$ can be canonically defined as the following chain of compositions of maps:
\begin{equation}\label{Eq:EChain}
\begin{aligned}
& \begin{tikzcd}
	\Aa \otimes \Bb_\omega & \Aa \otimes \Aa_R & \Aa^{\otimes \ZM} \otimes \Aa^{\otimes \ZM}  
		\arrow[rightarrow,from=1-1, to=1-2, "{\rm lift}"]
		\arrow[right ,from=1-2, to=1-3,"{\rm \tiny emb.}"]
\end{tikzcd}\\
& \qquad \qquad \qquad \qquad \qquad \qquad \begin{tikzcd}
	\ & \Aa^{\otimes \ZM} & \Aa_R \subset \Aa^{\otimes \ZM} &  \Bb_\omega,  
				\arrow[rightarrow ,from=1-1, to=1-2,"{\rm mult.}"]
	\arrow[rightarrow ,from=1-2, to=1-3,"{\rm shift}"]
	\arrow[rightarrow ,from=1-3, to=1-4,"{\rm proj.}"]
\end{tikzcd}
\end{aligned}
\end{equation}
where $\Aa$ is embedded in $\Aa_L \subset \Aa^{\otimes \ZM}$ in the most right position. Provided that we can make a rigorous sense of the above for generic nuclear algebras $\Aa$, we can try to apply \eqref{Eq:Facto0} and see if the outcome still reproduces the original state $\omega$. If we can provide sufficient conditions in which this happens, then we can enlarge the class of examples where states over tensor powers $\Aa^{\otimes \ZM}$ can be generated via rudimentary algorithms from the much simpler set of data $(\Aa,\Bb_\omega,\EM,\rho,e)$. 

\begin{remark}{\rm The authors of \cite{FannesCMP1992} showed us that any state $\omega$ can be reproduced by such an algorithm if we do not insist on the minimality of $\Bb$. In this work, we do insist that the data $(\Aa^{\otimes \ZM}, \omega)$ is fully reduced in the sense described above.
}$\Diamond$
\end{remark} 

When $\Aa$ is not finite dimensional, $\Bb_\omega$ can be infinite dimensional and this can happen even in some of the simplest scenarios. Indeed, let $\omega$ be a $C^\ast$-finitely correlated state over $\Aa^{\otimes \ZM}$ with reduced space $\Bb_\omega$ and maps $\EM$ and $\rho$. Now, let $\ZM'$ be a copy of $\ZM$ and consider $(\Aa^{\otimes \ZM'})^{\otimes \ZM}$ with the state supplied by 
\begin{equation}
\begin{tikzcd}
(\Aa^{\otimes \ZM'})^{\otimes \ZM} & (\Aa^{\otimes \ZM})^{\otimes \ZM'} & \CM.
\arrow[rightarrow, from=1-1, to=1-2,"\sim"]
\arrow[rightarrow, from=1-2, to=1-3,"\omega^{\otimes \ZM'}"]
\end{tikzcd}
\end{equation}
Hence, this example is about a vertical stacking of 1-dimensional spin chains. In this case, the reduced space is $\Bb_\omega^{\otimes \ZM'}$ and the map is $\EM^{\otimes \ZM'}$, while $\rho$ is amplified to $\rho^{\otimes \ZM'}$. This example is special because the assumption on $\Bb_\omega$ of being a finite $C^\ast$-algebra enables us to make sense of its infinite tensor power $\Bb_\omega^{\otimes \ZM'}$ as an AF-algebra. But we cannot expect this to happen in general cases. Therefore, our first outstanding task is to understand the structure of the reduced data and we call this phase of our program the reduction process. As we shall see in section~\ref{Sec:Redo}, the entanglement kernel is a kernel order ideal of $\Aa_R$ and, as a consequence, the quotient space $\Bb_\omega = \Aa_R/\Kk_\omega$ is a matrix order space with an order unit, which can be canonically Archimedeanized \cite{KavrukAM2013}. In other words, $\Bb_\omega$ always inherits an operator system structure from $\Aa_R$. Therefore, it comes equipped with the Archimedean order unit $e = 1 + \Kk_\omega$ and,  furthermore, if $q$ is the quotient map $\Aa_R \twoheadrightarrow \Bb_\omega$, then there exists a unique completely positive and unital map $\bar \omega$ such that $\omega = \bar \omega \circ q$. This map replaces $\rho$.

The second task in our program is making sense of the sequence~\eqref{Eq:EChain} and characterizing the resulted map $\EM$. As we have already seen, if the whole program is successful, then the original state factorizes through $\EM$ and, for this reason, we call this phase of our program the factorization process. Key to its progress are the results from \cite{KavrukJFA2011} on tensor products of operator systems. Specifically, if $\Aa$ is nuclear, Corollary~6.8 in \cite{KavrukJFA2011} assures us that the first tensor product seen in~\eqref{Eq:EChain} is unique but the operator system structure on the tensor product can be specified in many equivalent ways. We use the maximal tensor structure introduced in \cite{KavrukJFA2011} to prove that $\EM$ is a unital completely positive map.

The third task in our program is the investigation of a reconstruction algorithm based on \eqref{Eq:Facto0}, from an input data that shares similar properties with the reduced data of an actual state. In the context of infinite algebras, the algorithm produces products of operators, hence the name operator product states. While we show in section~\ref{Sec:Reconstruction} that such algorithm always produces a shift-invariant state on $\Aa^{\otimes \ZM}$ from such input data, the so obtained state may reduce to a different set of data than the input. To avoid such scenarios, we need a criterion to tell when the input data is actually the reduced data of some state, and such criterion is supplied in Theorems~\ref{Th:Main2} and \ref{Th:Main3}. We also supply sufficient conditions that ensure that the state is ergodic.

New examples of ergodic states on tensor powers $\Aa^{\otimes \ZM}$ can be constructed using just the technology developed in the already mentioned sections. In subsection~\ref{Sec:Reco}, the reader will find such example for the case when $\Aa$ is the group $C^\ast$-algebra of an infinite discrete amenable group. The use of Stinespring representations of the $\EM$ maps supplies additional routs to produce non-trivial examples of reconstructed states (see section~\ref{Sec:Stinespring}).

\section{Entanglement Kernel and the Reduced Space and State}

Throughout our presentation, we will oscillate between the categories of operator spaces and of operator systems, which are both extremely relevant for the program stated in our introduction. In this section, we first provide the minimal background needed to introduce the main concepts, formulate goals and sketch the road ahead. These initial steps can be formalized entirely in the category of operator spaces, hence we compiled a background material on it, mostly taken from the textbooks by Effros and Ruan \cite{EffrosBook}, by Blecher and Merdy \cite{BlecherBook} and by Pisier \cite{PisierBook}. It contains relevant definitions and fundamental statements that will be referenced throughout our presentation. This will will make the exposition self-sufficient and will fix the concepts and notation.

In the second part, we introduce and exemplify the main objects to be studied, namely, the algebra of physical observables, which is the tensor power $\Aa^{\otimes \ZM}$ with a nuclear $C^\ast$-algebra $\Aa$, the entanglement kernel of a given state $\omega$ over this algebra and the quotient space $\Bb_\omega$ of $\Aa^{\otimes \ZM}$ by this kernel. As we shall see, the latter has the structure of an operator space and $\omega$ descends to a completely contractive functional $\bar \omega$ over $\Bb_\omega$. We point out, however, that $\Bb_\omega$ has the potential to carry additional structure which anticipates the next steps for moving the program forward.

Before we start, let us lay out our conventions for the notation. The letters $H,K,L,\ldots$, will be designated for Hilbert spaces. The symbol $H^{(n)}$ will stand for the direct sum of $n$ identical copies of $H$, $H^{(n)}=H \oplus \ldots \oplus H$. The $C^\ast$-algebra of bounded linear maps between two Hilbert spaces will be denoted by $B(H,K)$ and, if $H$ coincides with $K$, the notation will be simplified to $B(H)$. The letters $E$, $F$, $G$, etc., will be designated to operator spaces. The matrix amplification of a linear space will be denoted by $M_n(E)$ and its elements will, most of the time, be indicated as $[e^{ij}]$ or $[e_{ij}]$. $M_{n,m}(\CM)$ will denote the space of linear maps from $\CM^m$ to $\CM^n$, equipped with the standard norm. The elements of the operator spaces and algebras will be denote by lowercase letters $e$, $f$, $g$, etc.. 

\subsection{Algebra of physical observables}
\label{Sec:PhysAlg}

Let $\Aa_j$, $j\in \ZM$, be $C^\ast$-algebras canonically isomorphic to a fixed unital and separable $C^\ast$-algebra $\Aa$, referred to as the site algebra. We denote by $\alpha_j: \Aa \twoheadrightarrowtail \Aa_j$ the canonical isomorphism. To avoid unnecessary complications, we assume that $\Aa$ is nuclear, such that all many possible ways to complete its algebraic tensor powers coincide. We introduce the notation
\begin{equation}
\Aa_{(m,n)} =\Aa_{m}\otimes \Aa_{m+1} \otimes \ldots \otimes \Aa_n, \quad m < n \in \ZM,
\end{equation}
and we will use the symbols $a_{(n,m)}$, $1_{(m,n)}$ and ${\rm id}_{(m,n)}$ for the generic monomials, the identity element and identity automorphism of $\Aa_{(m,n)}$, respectively. The natural unital embeddings 
\begin{equation}
\Aa_{(-n,n)} \rightarrowtail \Aa_{(-n-1,n+1)}, \quad a_{(-n,n)} \mapsto 1 \otimes a_{(-n,n)} \otimes 1,
\end{equation} 
supply an inductive tower
\begin{equation}
\Aa_0 \rightarrowtail \Aa_{(-1,1)} \rightarrowtail \cdots \rightarrowtail \Aa_{(-n,n)} \rightarrowtail \cdots
\end{equation}
of unital $C^\ast$-algebras, whose direct limit is the unital separable $C^\ast$-algebra $\Aa^{\otimes \ZM}$, denote here by $\Aa_\ZM$. This is the algebra of physical observables we are assuming in this work. It comes with canonical embeddings $\mathfrak i_{(m,n)} : \Aa_{(m,n)} \rightarrowtail \Aa_\ZM$.

Embedded in $\Aa_\ZM$, are the $C^\ast$ algebras $\Aa_{(n,\infty)}$ and $\Aa_{(-\infty,n)}$ defined by the inductive towers of $\Aa_{(n,m)}$ and $\Aa_{(m,n)}$ algebras, $m \to \pm \infty$, respectively. Special symbols will be used for $\Aa_R := \Aa_{(1,\infty)}$ and $\Aa_L : = \Aa_{(-\infty,0)}$. We note that $\Aa_\ZM = \Aa_L \otimes \Aa_R$ and also $\Aa_\ZM = \Aa_L \cdot \Aa_R$ when the latter are embedded in $\Aa_\ZM$, as well as that $\Aa_{R}$ ($\Aa_{L}$) belongs to the relative commutant of $\Aa_{L}$ ($\Aa_{R}$) inside $\Aa_\ZM$. 

As is the case for any $C^\ast$-algebra, $\Aa_\ZM$ comes equipped with a $C^\ast$-norm that enjoys the special property $\|a^\ast a\| = \|a^\ast \|\, \|a\| = \|a\|^2$, for any $a \in \Aa_\ZM$. Among many other things, this property enables one to define a special positive cone
\begin{equation}
\Aa^+_\ZM = \{ a^\ast a, \ a \in \Aa_\ZM\},
\end{equation} 
whose order semi-norm (see \ref{Def:ONorm}) is a (complete) norm and coincides with the $C^\ast$-norm. The state space of $\Aa_{\ZM}$ consists of all bounded linear functionals $\omega$ which map $\Aa^+_\ZM$ into $\RM_+$, the positive cone of $\CM$, and are normalized as $\omega(1)=1$. We will denote by $\Aa_R^+$ and $\Aa_L^+$ the positive cones of the corresponding $C^\ast$-algebras.

The algebra $\Aa_\ZM$ has a special (outer) automorphism $S:\Aa_\ZM \twoheadrightarrowtail \Aa_\ZM$, which is the shift acting on monomials as
\begin{equation}\label{Eq:ShiftOp}
S(\otimes a_n) = \otimes_{n \in \ZM} \, (\alpha_n \circ (\alpha_{n-1}^{-1})(a_{n-1}).
\end{equation}
Since $S$ shifts the entries from left to right, it maps $\Aa_R$ into itself, hence we can define a shift map $S_R$ on $\Aa_R$. A similar $C^\ast$-algebra morphism $S_L^{-1}$ can be defined on $\Aa_L$. The goal of our work is to explore the states $\omega$ on $\Aa_\ZM$ that are shift invariant, $\omega = \omega \circ S$, using the strategy develop in \cite{FannesCMP1992}.

\subsection{Background: Concrete and abstract operator spaces} Many classes of subspaces of $B(H)$ can be characterized concretely and abstractly, and operator spaces are no exceptions.

\begin{definition}[\cite{BlecherBook},~p.~5]\label{Def:OS} A concrete (closed) operator space is a (closed) linear subspace of $B(H)$ for some Hilbert space $H$.
\end{definition}

\begin{remark}{\rm The attribute ``closed'' is sometimes included in the definition of a concrete operator space, such as in \cite{PisierBook}, and sometimes is not, such as in \cite{EffrosBook} and \cite{BlecherBook}. We chose to side here with the second option and to specify explicitly when the encountered operator spaces are closed. Note that, in such cases, the operator spaces are actually complete.
}$\Diamond$
\end{remark}

Operator spaces are intrinsic structures that are intimately related to matrix amplifications of normed linear spaces:

\begin{definition}[\cite{EffrosBook},~p.~20] An abstract (closed) operator space is a linear space $E$ equipped with a system of (complete) matrix norms $\| \ \|_n$ on each $M_n(E)$, $n\in \NM^\times$ such that:
\begin{enumerate}[{\rm \ \ R1)}]

\item For all $e \in M_m(E)$ and $e' \in M_n(E)$,
\begin{equation}
\left \| \begin{pmatrix} e & 0 \\ 0 & e' \end{pmatrix} \right \|_{m+n} = {\max}\big \{\|e\|_n, \, \| e'\|_m\big \};
\end{equation}

\item For all $e \in M_m(E)$, $\alpha \in M_{n,m}(\CM)$ and $\beta \in M_{m,n}(\CM)$,
\begin{equation}
\|\alpha e \beta\|_n \leq \|\alpha\| \|e\|_m \|\beta\|. 
\end{equation}

\end{enumerate}
\end{definition}

\begin{theorem}[\cite{RuanJFA1988,EffrosPAMS1993}]\label{Th:Ruan1} Any abstract operator space can be isometrically embedded in the $B(H)$ of some Hilbert space $H$. Conversely, if $E$ can be isometrically embedded in $B(H)$, then the norms $\| \ \|_n$ inherited by $M_n(E)$ from $M_n(B(H)) \simeq B(H^{(n)},H^{(n)})$ satisfy {\rm R1} and {\rm R2}.
\end{theorem} 

It is important to expose the fine synergies set in motion by the conditions R1 and R2, as revealed by the following fact:

\begin{proposition}[\cite{EffrosBook},~p.~22 and 34]\label{Pro:RuanVsNorms} Suppose that $E$ is a linear space, and that we are provided with \underline{mappings} $\| \cdot \|_n : M_n(E) \to [0,\infty)$
for all $n \in \NM^\times$, satisfying {\rm R1} (or a slightly weaker version) and {\rm R2}. Then these mappings are seminorms which satisfy {\rm R1} and {\rm R2}. If, in addition, $\| \cdot \|_1$ is a (complete) norm, then the same is true for all matrix seminorms.
\end{proposition}

Closed linear sub-spaces of an operator space are again operator spaces with the $n$-norms induced from the parent operator space. More importantly for us is a fundamental result by Ruan that quotients of operator spaces by closed linear sub-spaces are also operators spaces:

\begin{proposition}[\cite{BlecherBook}~p.~8, \cite{EffrosBook}~Prop.~3.1.1]\label{Pro:Quotient1} If $E$ is an operator space and $F$ is one of its closed linear subspaces, then $E/F$ is an operator space with norms induced by the identification $M_n(E/F) \simeq M_n(E)/M_n(F)$. Explicitly, these norms are give by the formula
\begin{equation}
\|[e_{ij}+ F] \|_n =\inf \big \{\| [e_{ij}+f_{ij}]\|_n, \ [f_{ij}] \in M_n(F)\big \} , 
\end{equation}
for any $[e_{ij}] \in M_n(E)$.
\end{proposition} 

\begin{remark}\label{Re:Closed}{\rm Clearly, if $E$ is closed, then $\|\cdot \|_1$ defined above is complete. Then Proposition~\ref{Pro:RuanVsNorms} assures us that all $\| \cdot \|_n$ norms are complete.
}$\Diamond$
\end{remark}

Concrete presentations of quotient operator spaces were supplied by Rieffel in \cite{RieffelMS2014}. Unfortunately, we were not able to take advantage of them at this point.

\subsection{Background: Completely bounded linear maps}

The morphisms in the category of operator spaces are supplied by the completely bounded (c.b.) linear maps:

\begin{definition}[\cite{PisierBook}~p.~19; \cite{BlecherBook}~p.~4] A linear map $u:E \rightarrow F$ between two operator spaces can be amplified to a linear map
\begin{equation}
u_n : M_n(E) \rightarrow M_n(F), \quad u_n([e_{ij}]) =[u(e_{ij})],
\end{equation}
 for all $n \geq1$. The map $u$ is called:
\begin{enumerate}[{\rm \ \ 1)}]
\item Completely bounded if
\begin{equation}
\sup_{n \geq 1} \, \|u_n\|_{M_n(E) \rightarrow M_n(F)}  < \infty.
\end{equation}
\item Complete isometry if all $u_n$'s are isometries.
\item Complete quotient if each $u_n$ sends the unit ball of $M_n(E)$ onto the unit ball of $M_n(F)$.
\end{enumerate}
\end{definition}

The set of c.b. maps ${\rm CB}(E,F)$ is closed under addition and becomes a Banach linear space when equipped with the norm
 \begin{equation}
\|u\|_{\rm cb}=\sup_{n \in \NM} \, \|u_n\|_{M_n(E) \rightarrow M_n(F)}.
\end{equation}
As expected, c.b. linear maps behave well under composition:

\begin{proposition}[\cite{PisierBook}~p.~19]\label{Pro:CBComp} If $E$, $F$ and $G$ are operator spaces and $u:E \rightarrow F$ and $v: F \rightarrow G$ are completely bounded linear maps, then $v \circ u : E \rightarrow G$ is a completely bounded map and $\|v \circ u\|_{\rm cb} \leq \|v\|_{\rm cb} \, \|u\|_{\rm cb}$.
\end{proposition}

C.b. linear maps also behave quite natural under taking quotients:

\begin{proposition}[\cite{PisierBook}~p.~42]\label{Pro:Quotient2} Let $E$, $F$ and $G$ be operator spaces such that $F \subset E$, and let $q: E \rightarrow E/F$ be the canonical surjection. Then, a linear map $u:E/F \rightarrow G$ is completely bounded if and only if $u\circ q$ is completely bounded and $\|u\|_{\rm cb} = \|u\circ q\|_{\rm cb}$.
\end{proposition}

\begin{proposition}[\cite{BlecherBook}~p.~8]\label{Pro:Quotient3} If $u:E \rightarrow G$ is completely bounded and if $F$ is a closed subspace of $E$ contained in ${\rm Ker} \, u$, then the canonical map $\tilde u : E/F \rightarrow G$ induced by $u$ is also completely bounded, with $\| \tilde u \|_{\rm c.b.} = \| u \|_{\rm c.b.}$. If $F = {\rm Ker} \, u$, then $u$ is a complete quotient map if and only if $\tilde u$ is a complete isometric isomorphism.
\end{proposition}

\begin{corollary}\label{Cor:Q} Let $E$ and $F$ be operator spaces with $F \subset E$. Then the canonical surjection $q:E \rightarrow E/F$ is a complete quotient map and ${\rm Ker} \, 	q = F$.
\end{corollary}

The following statement is known as the fundamental factorization/extension of c.b. linear maps.

\begin{theorem}[\cite{PisierBook}~p.~23]\label{Th:CB} 
Consider a completely bounded map
\begin{equation}
\begin{tikzcd}
	B(H) & B(K) \\
	E   & F
		\arrow[rightarrowtail,from=2-1, to=1-1]
		\arrow[rightarrow ,from=2-1, to=2-2,"u"]
	\arrow[rightarrowtail ,from=2-2, to=1-2]
\end{tikzcd}
\end{equation}
Then there exits a Hilbert space $\widehat H$, a representation $\pi : B(H) \rightarrow B(\widehat H)$ and operators $V_1 : K \rightarrow \widehat H$, $V_2 : \widehat H \rightarrow K$ such that $\|V_1\|\, \|V_2\| = \|u\|_{\rm cb}$ and
\begin{equation}\label{Eq:Facto1}
u(e) = V_2 \pi(e) V_1, \quad \forall \ e\in E.
\end{equation}
Conversely, if \eqref{Eq:Facto1} holds, then $u$ is completely bounded and $\|u\|_{\rm cb} \leq \|V_1\|\, \|V_2\|$. Moreover, $u$ admits a completely bounded extension $\tilde u$ such that
\begin{equation}
\begin{tikzcd}
	B(H) & B(K) \\
	E   & F
		\arrow[rightarrowtail,from=2-1, to=1-1]
		\arrow[rightarrow ,from=2-1, to=2-2,"u"]
		\arrow[rightarrow ,from=1-1, to=1-2,"\tilde u"]
	\arrow[rightarrowtail ,from=2-2, to=1-2]
\end{tikzcd}
\end{equation}
is a commutative diagram and $\|\tilde u\|_{\rm cb} = \|u\|_{\rm cb}$.
\end{theorem}

\subsection{The entanglement kernel and the reduced state}

Let $\omega$ be a state on $\Aa_\ZM$, not necessarily shift invariant. Then, for any $x \in \Aa_L \subset \Aa_Z$, define a bounded linear functional
\begin{equation}
\omega_x : \Aa_R \subset \Aa_\ZM \rightarrow \CM, \quad \omega_x(a_R) = \omega(x \, a_R),
\end{equation} 
which is not positive, in general. Inspired by \cite{FannesCMP1992}, we introduce:

\begin{definition}\label{Def:EntKer} The following subset of $\Aa_R$,
\begin{equation}\label{Eq:Kappa1}
\Kk_\omega = \bigcap_{x \in \Aa_L} {\rm Ker} \, \omega_x,
\end{equation}
will be referred to as the entanglement kernel of $\omega$. 
\end{definition}

$\Kk_\omega$ is an intersection of closed linear sub-space, hence it is closed linear sub-space of the $C^\ast$-algebra $\Aa_R$. As such, it is automatically a closed operator subspace and it enters the exact sequence of closed operator spaces
\begin{equation}\label{Eq:Seq}
\Kk_\omega \rightarrowtail  \Aa_R 
 \twoheadrightarrow \Bb_\omega = \Aa_R / \Kk_\omega \ .
\end{equation}
Indeed, from Proposition~\ref{Pro:Quotient1} and Remark~\ref{Re:Closed}, we know that the quotient space $\Bb_\omega$ inherits a natural closed operator space structure. Furthermore, the second map in \eqref{Eq:Seq} is the canonical surjection $q: \Aa_R \to \Aa_R/\Kk_\omega$, which is a complete quotient map, as we learned from Corollary~\ref{Cor:Q}. We will refer to $\Bb_\omega$ as the $\omega$-reduced operator space, which can be entirely and abstractly described by the data $\big (\Bb_\omega,\{\|\cdot \|^{\rm osp}_n\}_{n\geq 1}\big )$, with the matrix norms supplied by Proposition~\ref{Pro:Quotient1}. Its elements will be specified by $b$, $b'$ and so on. Also, the matrix amplifications of $q$, which are all contractions, will be denoted by $q_n$. The class of an element $a_R \in \Aa_R$ in $\Bb_\omega = \Aa_R / \Kk_\omega$ will be indicated by several symbols, such as
\begin{equation}
q(a_R) = \hat a_R = \lfloor a_R \rfloor.
\end{equation}
The second notation $\hat a_R$ is useful when considering matrix amplifications of $\Bb_\omega$. The third notation will be used when $a_R$ is given as a long expression. 

\begin{proposition}\label{Pro:Omega} Let $\omega_R$ be the state on $\Aa_R$ supplied by the restriction of $\omega$. Then $\omega_R$ descends to a completely bounded linear functional $\bar \omega : \Bb_\omega \rightarrow \CM$ with $\omega_R = \bar \omega \circ q$ and $\|\bar \omega\|_{\rm cb}=1$.
\end{proposition}

\proof Taking $x$ the unit of $\Aa_L$ sub-algebra, we see from the definition \eqref{Eq:Kappa1} of $\Kk_\omega$ that $\Kk_\omega \subset {\rm Ker} \, \omega_R$. As such, the map
\begin{equation}
\bar \omega : \Bb_\omega \rightarrow \CM, \quad \bar \omega(a_R + \Kk_\omega) = \omega_R(a_R), 
\end{equation}
is well defined. Furthermore, $\bar \omega \circ q = \omega_R$ and the latter is a completely bounded functional with c.b. norm 1. According to Proposition~\ref{Pro:Quotient2}, this can be true if and only if $\bar \omega$ is completely contractive. \qed

\subsection{Examples of entanglement kernels and reduced spaces}

The algebra $\Aa_R$ can be reduced (quotiented) in many different ways, but one of the practical values of the above particular reduction, which is the great insight supplied by \cite{FannesCMP1992}, is that $\Kk_\omega$ and $\Bb_\omega$ can be computed for a large class of interesting physical models.

\begin{example}\label{Ex:ProdState1}{\rm Let $\omega_0$ be a state on $\Aa$ and let $\omega=\omega_0^{\otimes \ZM}$ be the product state on $\Aa_\ZM$. In this case, 
\begin{equation}
\omega_x(a_R) = \omega(x)\, \omega(a_R) = \omega(x) \, \omega_R(a_R)
\end{equation} 
for any $x \in \Aa_L$ and $a_R \in \Aa_R$, hence, $\Kk_\omega = {\rm Ker}\, \omega_R$. Since ${\rm Ker}\, \omega_R$ coincides with the linear subspace $\{a_R -\omega_R(a_R) \cdot 1,  \ a_R \in \Aa_r\}$, it follows that $\Bb_\omega = \CM \cdot 1$ for any product state on $\Aa_\ZM$.
}$\Diamond$
\end{example}

Any product state has zero correlation length, in the language introduced in \cite{FannesCMP1992}. Of course, the main interest of the physics community is on correlated states. The following example supplies a large class of such states for which $\Bb_\omega$ is again computable, at least formally.

\begin{figure}[t]
\includegraphics[width=0.5\linewidth]{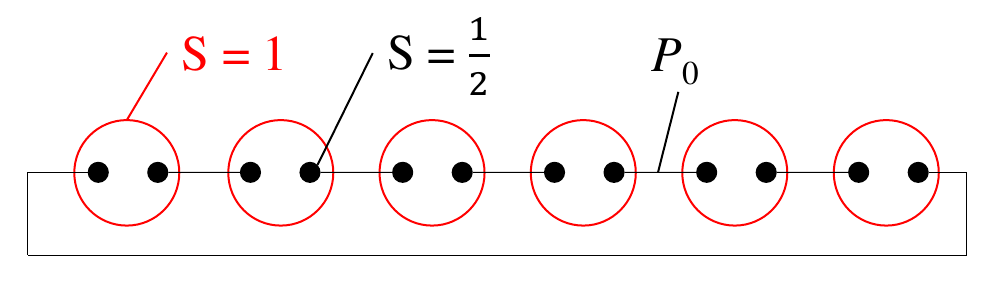}
\caption{\small Spin-1 particles, shown by red bubbles, are arranged in a closed chain of length $N$. The algebra of observables is $[1]^{\otimes \ZM_N}$, where $[j]$ denotes the spin-j algebra. The site algebra $[1]$ is embedded via a map $\mathfrak j$ into the algebra $\big [ \tfrac{1}{2} \big ] \otimes \big [\tfrac{1}{2} \big ]$ of two spin-$\frac{1}{2}$ particles, shown as pairs of black dots, and the algebra $[1]^{\otimes \ZM_N}$ is embedded into $[\frac{1}{2}]^{\otimes \ZM_{2N}}$ via $\mathfrak j^{\otimes \ZM_N}$. The tensor power $P_0^{\otimes \ZM_N}$ of the projection $P_0$ onto the one-dimensional subspace of the decomposition $\big [ \tfrac{1}{2} \big ] \otimes \big [ \tfrac{1}{2} \big ] \simeq [0] \oplus [1]$, shown by the segments, generates a rank-one projection in $[1]^{\otimes \ZM_N}$ and a state $[1]^{\otimes \ZM_N} \ni M \mapsto {\rm Tr}\big ( \mathfrak j^{\otimes \ZM_N}(M)P_0^{\otimes \ZM_N} \big )$. Note that the projections are applied on the ``bonds'' shown by the sticks and this is why the shift appears in \eqref{Eq:AKLTSeq1}. }
\label{Fig:1DSystem}
\end{figure}

\begin{example}\label{Ex:AKLT1}{\rm Product states in conjunction with shift maps can be used in  creative ways to generate states with finite correlation length. The following class of states is modeled after the so call AKLT state for spin-1 system  \cite{AffleckCMP1988}, whose construction is sketched in Fig.~\ref{Fig:1DSystem}. In fact, the construction given here covers all dimerized states introduced in \cite{AffleckCMP1988}[p.~523]. For this reason, we refer to the states cover by this example  as AKLT-type. The construction involves a (nuclear) $C^\ast$-algebra $\tilde \Aa$ and a projection $p \in  \tilde \Aa \otimes \Aa$. Then the local algebra $\Aa$ is defined as the unital $C^\ast$-algebra $\Aa = p (\tilde \Aa \otimes \tilde \Aa) p$, with $p$ playing the role of the unit. We generate a state $\omega$ on $\Aa_\ZM$ via the thermodynamic limit $N \rightarrow \infty$ of the states $\omega_N$ on $\Aa_{\ZM_N}$, supplied by the following sequence of maps
\begin{equation}\label{Eq:AKLTSeq1}
\begin{tikzcd}
	\Aa_{\ZM_N} & (\tilde  \Aa \otimes \tilde \Aa)_{\ZM_N} \simeq  \tilde \Aa_{\ZM_{2N}} &  \tilde \Aa_{\ZM_{2N}} \simeq \big (\tilde \Aa \otimes \tilde \Aa \big )_{\ZM_N} & \CM,
		\arrow[rightarrowtail,from=1-1, to=1-2,"\mathfrak j^{\otimes \ZM_N}"]
		\arrow[rightarrow ,from=1-2, to=1-3,"S"]
	\arrow[rightarrow ,from=1-3, to=1-4,"\xi_0^{\otimes \ZM_N}"]
\end{tikzcd}
\end{equation}
where $\xi_0$ is a positive functional on $\tilde \Aa \otimes \tilde \Aa$. The first map is the power of the {\it non}-unital inclusion
\begin{equation}\label{Eq:NUIncl}
\Aa \ni p \hat a p \mapsto \mathfrak j (p \hat a p) : = p \hat a p \in \tilde \Aa \otimes \tilde \Aa,
\end{equation} 
and $S$ is the obvious shift map on $\tilde \Aa_{\ZM_{2N}}$. 

If Eq.~\eqref{Eq:AKLTSeq1} is to supply a state on $\Aa_{\ZM}$, the above data must obey the constraint
\begin{equation}
\lim_{N \rightarrow \infty} \xi_0^{\otimes \ZM_N} \Big ( S \big (p^{\otimes \ZM_N} \big ) \Big ) = 1.
\end{equation}
For example, this is indeed satisfied if $({\rm id} \otimes \xi_0)(p \otimes \tilde 1) = \tilde 1$ and, in particular, for the case described in Fig.~\ref{Fig:1DSystem}. The limiting procedure is required because $p^{\otimes \ZM}$ is not an element of $(\tilde \Aa \otimes \tilde \Aa)^{\otimes \ZM}$, hence the non-unital inclusion $\mathfrak j^{\otimes \ZM_N}$ does not make sense in the thermodynamic limit. Let us specify that, if a unital inclusion is used instead of $\mathfrak j$ and $\xi_0$ is a state, then all the maps in the sequence \eqref{Eq:AKLTSeq1} are well defined in the thermodynamic limit, but then the state on $\Aa$ is a trivial product state. We will not investigate the thermodynamic limiting process here because we will re-construct this class of states  via a different path in Example~\ref{Ex:AKLT2}. We only mention here that these  issues  were fully resolved in \cite{AffleckCMP1988} for the particular case illustrated in Fig.~\ref{Fig:1DSystem}. 

Now, assuming that the thermodynamic limit of the state exists, we compute the corresponding operator spaces $\Kk_\omega$ and $\Bb_\omega$. For this, we first note the obvious isomorphisms
\begin{equation}
\chi_R : (\tilde \Aa \otimes \tilde \Aa)_R \twoheadrightarrowtail \tilde \Aa \otimes (\tilde \Aa \otimes \tilde \Aa)_R, \quad \chi_L : (\tilde \Aa \otimes \tilde \Aa)_L \twoheadrightarrowtail (\tilde \Aa \otimes \tilde \Aa)_L \otimes \tilde \Aa,
\end{equation}
which we use to define the maps
\begin{equation}
\Gamma_R = ({\rm id} \otimes \xi_0^{\otimes \NM^\times}) \circ \chi_R \circ \mathfrak j^{\otimes \NM^\times}, \quad \Gamma_L = (\xi_0^{\otimes(\ZM \setminus \NM^\times)} \otimes {\rm id}) \circ \chi_L \circ \mathfrak j^{\otimes (\ZM\setminus \NM^\times)},
\end{equation}
from $\Aa_R$ to $\tilde \Aa$ and from $\Aa_L$ to $\tilde \Aa$, respectively. These maps are well defined because of our assumption that the map $\xi_0^{\otimes \ZM} \circ S \circ \mathfrak j^{\otimes \ZM}$ exists as the thermodynamic limit of the chain of maps~\eqref{Eq:AKLTSeq1}. Also, there exists a specific isomorphism between $\Aa_L$ and $\Aa_R$, induced by the reflection of $\ZM$ relative to $\frac{1}{2}$, which sends $\Gamma_{L,R}$ into $\Gamma_{R,L}$, respectively. From this fact, we deduce that the ranges of $\Gamma_L$ and $\Gamma_R$ coincide. Furthermore, we have
 \begin{equation}
\omega(a_L a_R) = \xi_0(\tilde a_L \otimes \tilde a_R), \quad \tilde a_{R,L} = \Gamma_{R,L}(a_{R,L}).
\end{equation}
One immediate conclusion from here is that ${\rm Ker}\, \Gamma_R \subseteq \Kk_\omega$. If, additionally, we choose $\xi_0$ such that $\xi_0(\tilde a \otimes \tilde a) \neq 0$ when $\tilde a$ samples a dense sub-set of $\tilde \Aa$, then ${\rm Ker}\, \Gamma_R = \Kk_\omega$, because, for each $\tilde a_R$, we can produce an $\tilde a_L$ via $\Gamma_L$ that coincide with $\tilde a_R$. The conclusion is that $\Bb_\omega = {\rm Range}\, \Gamma_R=\tilde \Aa$ for such choice of $\xi_0$. 
}$\Diamond$
\end{example}

\subsection{Additional structures and a look ahead}
\label{Sec:RS}

Let us denote by $\bar \omega_n$ the matrix amplifications of the reduced state $\bar \omega$. We want to point here to a few interesting properties of these maps. Specifically, since $M_n(\Aa_R)$ are again $C^\ast$-algebras, we can consider their positive cones and define the subspaces
\begin{equation}\label{Eq:Dn}
\Dd_n : = M_n(\Aa_R)^+/M_n(\Kk_\omega) = q_n\big (M_n(\Aa_R)^+\big ),
\end{equation} 
which can be characterized more explicitly as \cite{KavrukAM2013}[p.~327]
\begin{equation}
\Dd_n = \big \{ [a_{ij}+\Kk_\omega] \in M_n(\Bb_\omega) \ | \ \exists \, k_{ij} \in \Kk_\omega \ s.t. \ [a_{ij}+k_{ij}] \in M_n(\Aa_R)^+ \big \}.
\end{equation}
Evidently, we have
\begin{equation}\label{Eq:BarOm}
\bar \omega_n (\Dd_n) \subset \RM_+, \quad n=1,2, \ldots .
\end{equation}

This indicates that $\Bb_n$ spaces may carry more structure, perhaps that of an operator system, but, without a more refined characterization of the entanglement kernel, this cannot be established. Let us point out that, as of now, the structure of the reduced data is quite far from the one assumed in \cite{FannesCMP1992}, where the authors focused on the special cases where $\Bb_\omega$ is a $C^\ast$-algebra and the reduced state $\bar \omega$ is a completely positive map. As we shall see, in general, $\Dd_n$'s defined above do not generate an Archimedean matrix order structure, hence an abstract operator structure. Here are a few things that can go wrong:
\begin{itemize}
\item $\Dd_n$'s may fail to be closed spaces; 
\item $\Dd_n$'s may not be compatible, in the sense that the relation $A^\ast \Dd_n A \subseteq \Dd_m$ may fail for $A$ an ordinary $n \times m$ matrix;
\item The intersections $\Dd_n \cap (-\Dd_n)$ might contain elements other than $0$. 
\end{itemize}
Of course, there are states for which $\Dd_n$'s do supply Archimedean matrix order structures and explicit sufficient conditions for this to happen will be supplied in the next section. 

To summarize, the data that we pass to the reduction process is that of an operator space with extra structure:
\begin{equation}
\{\Bb_\omega,\{\|\cdot \|_n^{\rm osp}\}_{n \geq 1}, \{\Dd_n\}_{n\geq 1}, 1+\Kk_\omega,\bar \omega: \Bb_\omega \to \CM\}.
\end{equation}

\section{Reduction Process}
\label{Sec:Redo}

A $C^\ast$-algebra comes equipped with an order structure that is compatible with the topology induced by the $C^\ast$-norm (see Remark~\ref{Re:CStarOrder}). As a closed linear subspace of $C^\ast$-algebra $\Aa_R$, we have seen that $\Kk_\omega$ inherits a closed operator space structure and so does its corresponding quotient space $\Bb_\omega$. As already hinted, these spaces also inherit a matrix order structure that can be Archimedeanized to an operator system structure. The latter induces an order topology on $\Bb_\omega$ that in general is different from the operator space topology \cite{KavrukAM2013}[Sec.~4]. Therefore, a decision must be made about which of the two inherited structures is more important for the present context. Our main goal for this section is to expose and characterize the inherited order structures, which is achieved over the course of several subsections. An additional section will discuss the relation between the operator space and the operator system structures on $\Bb_\omega$ and will give sufficient conditions that assure that the two structures coincide. For example, this is always the case if the inherited operator system structure is close and, as such, its matrix order norms are complete. Lastly, a choice will be made in favor of the operator system structure on $\Bb_\omega$ and, with that, we can finally describe what we call the reduced data.

The proof of the existence of a canonical Archimedean matrix order structure on $\Bb_\omega$ consists of tying together several concepts and results from the existing literature, due to Kavruk, Paulsen, Todorov and Tomforde \cite{PaulsenIUMJ2009,PaulsenPLMS2010,KavrukJFA2011,KavrukAM2013}. We will take this opportunity and give a brisk recap of these ideas, which supply the natural framework and the right tools for the problem at hand, something that we still contemplate with amazement.\footnote{Specifically, that the pioneering concepts introduced in \cite{FannesCMP1992} found their rightful home in a framework developed two decades later.}  In the process, one will hear about ordered $\ast$-vector spaces, (Archimedean) order units, order semi-norms and topologies, as well as order ideals and their quotients \cite{PaulsenIUMJ2009}. One will also hear about matrix ordered $\ast$-vector spaces, (Archimedean) matrix order units and a conceptual refinement of the order ideal, which is the kernel introduced in \cite{KavrukAM2013}. The later has the remarkable property that its quotient space carries automatically an ordered $\ast$-vector space with an Archimedean matrix order unit.

\subsection{Background: Ordered vector spaces and their order topologies}\label{Sec:OrderSpaces} This material, which is entirely collected from \cite{PaulsenIUMJ2009}, will help us elucidate the structure of $\Bb_\omega$, as induced by the subsets $\Dd_n$ introduced in the previous section.

\begin{definition}[\cite{PaulsenIUMJ2009},~p.~1322] If $\Vv$ is a real vector space, a cone in $\Vv$ is a nonempty subset $\Cc \subseteq \Vv$ with the following two properties:
\begin{enumerate}[{\ \ \rm 1)}]
\item $a v \in \Cc$ whenever $a \in [0,\infty)$ and $v \in \Cc$;
\item $v + w \in \Cc$ whenever $v,w \in \Cc$.
\end{enumerate}
An ordered vector space is a pair $(\Vv,\Vv^+)$ consisting of a real vector space $\Vv$ and cone $\Vv^+ \subseteq \Vv$ satisfying $\Vv^+\cap (-\Vv^+) = \{0\}$.
\end{definition}

\begin{remark}{\rm If $(\Vv, \Vv^+)$ is an ordered real vectors space, one writes $v \geq v'$ if $v-v' \in \Vv^+$.
}$\Diamond$
\end{remark}

\begin{definition}[\cite{PaulsenIUMJ2009},~pp.~1323-1324] If $(\Vv,\Vv^+)$ is an ordered real vector space, an element $e \in \Vv$ is called an order unit for $\Vv$ if, for each $v \in \Vv$, there exists a real number $r >0$ such that $r e \geq v$. The order unit $e$ is called Archimedean if whenever $v \in \Vv$ with $r e + v \geq 0$ for all real $r>0$, then $v \in \Vv^+$.
\end{definition}

\begin{example}[\cite{PaulsenIUMJ2009},~p.~1353]\label{Ex:OrderAlg}{\rm The real vector space of self-adjoint elements of any unital $C^\ast$-algebra is an ordered vector space with the unit playing the role of Archimedean order unit.}\ $\Diamond$
\end{example}

Of course, our interest is in order structures on complex vector spaces. In this case, an extra structure is required.

\begin{definition}[\cite{PaulsenIUMJ2009},~p.~1337] A $\ast$-vector space consists of a complex vector space $\Vv$ together with a map $\ast : \Vv \rightarrow \Vv$ that is involutive, $(v^\ast)^\ast = v$ for all $v \in \Vv$, and conjugate linear. If $\Vv$ is a $\ast$-vector space, then $\Vv_{\rm h}=\{v \in \Vv \ | \ v^\ast = v\}$ represents the set of hermitean elements of $\Vv$. It carries the structure of a real vector space.
\end{definition}

\begin{definition}[\cite{PaulsenIUMJ2009},~p.~1337] If $\Vv$ is a $\ast$-vector space, one says that $(\Vv,\Vv^+)$ is an ordered $\ast$-vector space if $(\Vv_{\rm h},\Vv^+)$ is an ordered real vector space. Furthermore, $e \in \Vv$ is an (Archimedean) order unit for $(\Vv,\Vv^+)$ if it is an (Archimedean) order unit for $(\Vv_{\rm h},\Vv^+)$.
\end{definition}

\begin{definition}[\cite{PaulsenIUMJ2009},~p.~1337] Let $(\Vv,\Vv^+)$ be an ordered $\ast$-vector space with order unit $e$ and let $(\Ww,\Ww^+)$ be an ordered $\ast$-vector space with order unit $e'$. A linear map $\varphi: \Vv \rightarrow \Ww$ is positive if $v \in \Vv^+$ implies $\varphi(v) \in \Ww^+$, and unital if $\varphi(e)=e'$.
\end{definition}

Order structures can be used to generate topologies:

\begin{definition}[\cite{PaulsenIUMJ2009},~p.~1327]\label{Def:ONorm} Let $(\Vv,\Vv^+)$ be an ordered real vector space with order unit $e$. The order semi-norm on $\Vv$ determined by $e$ is defined as:
\begin{equation}
\llbracket v \rrbracket = \inf \{r \in \RM \, | \, re + v \geq 0 \ \mbox{and} \ re - v \geq 0 \}.
\end{equation}
The order topology on $\Vv$ is the topology induced by the order semi-norm.
\end{definition}

The following statement supplies a complete characterization of the order seminorm:

\begin{theorem}[\cite{PaulsenIUMJ2009},~p.~1330] Let $(\Vv,\Vv^+)$ be an ordered real vector space with order unit $e$. Then the order seminorm $\llbracket \cdot \rrbracket$ is the unique seminorm on $\Vv$ satisfying simultaneously the following three conditions:
\begin{enumerate}[{\rm \ \ 1)}]
\item $\llbracket e \rrbracket = 1$;
\item If $-v' \leq v \leq v'$, then $\llbracket v \rrbracket\leq \llbracket v' \rrbracket$;
\item If $f: \Vv \rightarrow \RM$ is a state, then $|f(v)| \leq \llbracket v \rrbracket$.
\end{enumerate}
\end{theorem}

When the order unit is Archimedean, then $\llbracket \cdot \rrbracket$ is actually a norm and the order topology is Hausdorff (the reciprocal is not necessarily true \cite{PaulsenIUMJ2009}[p.~1328]). Nevertheless, the Archimedean case can be characterized as it follows:

\begin{theorem}[\cite{PaulsenIUMJ2009},~p.~1330]\label{Th:Order} Let $(\Vv,\Vv^+)$ be an ordered real vector space with order unit $e$, and let $\llbracket \cdot \rrbracket$ be the order semi-norm determined by $e$. Then the following are equivalent:
\begin{enumerate}[{\ \ \rm i)}]
\item $e$ is Archimedean;
\item $\Vv^+$ is a closed subset of $\Vv$ in the order topology induced by $\llbracket \cdot \rrbracket$;
\item $-\llbracket v\rrbracket \, e \leq v \leq \llbracket v\rrbracket \, e$ for all $v \in \Vv$.
\end{enumerate}
\end{theorem}

\begin{remark}[\cite{PaulsenIUMJ2009},~Sec.~4]{\rm The order semi-norm on the hermitean sub-space of an ordered $\ast$-vector space with unit can be extend over the entire complex space, in an essentially unique way.
}$\Diamond$
\end{remark}

\begin{remark}\label{Re:CStarOrder}{\rm $C^\ast$-algebras are extremely special cases where the $C^\ast$-norm coincides with order norm. In particular, this is always the case for $B(H)$ of a Hilbert space.
}$\Diamond$
\end{remark} 
 
\begin{definition}[\cite{PaulsenIUMJ2009},~p.~1341]\label{Def:OrderIdeal} If $(\Vv,\Vv^+)$ is an ordered $\ast$-vector space, then a subspace $\Jj \subseteq \Vv$ is called an order ideal if $\Jj$ is self-adjoint ($\Jj^\ast=\Jj)$ and, furthermore, $v\in \Jj \cap \Vv^+$ and $0\leq v' \leq v$ implies that $v' \in \Jj$.
\end{definition}

\begin{proposition}[\cite{PaulsenIUMJ2009},~p.~1342]\label{Pro:OrderQ} Let $(\Vv,\Vv^+)$ be an ordered $\ast$-vector space with order unit $e$ and let $\Jj \subset V$ be an order ideal. Then $(\Vv/\Jj,\Vv^+ / \Jj)$ is an ordered $\ast$-vector space with order unit $e + \Jj$.
\end{proposition}

\subsection{Entanglement kernel is an order ideal} We prove this statement in several steps.

\begin{proposition} $\Kk_\omega$ does not contain the unit $1_R$ of $\Aa_R$.
\end{proposition}

\begin{proof} We need to find one element $x$ of $\Aa_L$ for which $\omega_x(1_R)=\omega(x) \neq 0$. This element is the identity of $\Aa_L$. \end{proof}

\begin{proposition} The entanglement kernel is self-adjoint: $\Kk_\omega^\ast = \Kk_\omega$.
\end{proposition}

\proof We have
\begin{equation}
\omega_x(a_R^\ast)=\omega(x \, a_R^\ast) = \omega( a_Rx^\ast)^\ast=\omega(x^\ast a_R)^\ast,
\end{equation}
for all $\ x \in \Aa_L$, where for the last equality we use that $[\Aa_L,\Aa_R]=\{0\}$. Hence, if $a_R \in \Kk_\omega$, then $\omega(x a_R^\ast)=0$ for any $x \in \Aa_L$. As a consequence, $a_R^\ast \in \Kk_\omega$.\qed

\begin{proposition}\label{Pro:K2} The entanglement kernel can be equivalently defined as
\begin{equation}\label{Eq:Kappa2}
\Kk_\omega = \bigcap_{x \in \Aa_L^+} {\rm Ker} \, \omega_x.
\end{equation}
Compared to \eqref{Eq:Kappa1}, the intersection in \eqref{Eq:Kappa2}  runs over the (much) smaller space of positive elements.
\end{proposition}

\proof Let $\Kk'_\omega$ denote the set at the right side of \eqref{Eq:Kappa2}. Clearly, $\Kk_\omega \subseteq \Kk'_\omega$. Consider now an $a_R$ from $\Kk'_\omega$ and let $x$ be an arbitrary element from $\Aa_L$. The latter can be always  decompose as $x=(x_1^+  - x_2^+) + \imath (x_3^+ - x_4^+)$ in terms of positive elements $x_i^+ \in \Aa_L^+$, though this decomposition is not unique. Nevertheless, from \eqref{Eq:Kappa2}, we know that $\omega(x_i^+ a_R)=0$, for all $i=\overline{1,4}$. Since $\omega$ is a linear map, we can conclude that $\omega(x\, a_R)=0$, hence $a_R$ also belong to the set defined at\eqref{Eq:Kappa1}.\qed

\begin{proposition}\label{Pro:KIdeal} The entanglement kernel is an order ideal of $\Aa_R$.
\end{proposition}

\proof Consider $k_R^+$ from $\Kk_\omega \cap \Aa_R^+$ and $a_R^+$ from $\Aa_R^+$, such that $a_R^+ \leq k_R^+$. Our task is to show that $a_R^+$ is automatically in $\Kk_\omega \cap \Aa_R^+$. For this, consider $x^+ \in \Aa_L^+$ and note that $x^+ a_R^+$ and $x^+ k_R^+$ are positive element of $\Aa_\ZM$, because $x^+$ commutes with any element from $\Aa_L$, and clearly $x^+ a_R^+ \leq x^+ k_R^+$. Then:
\begin{equation}
0 \leq \omega(x^+ a_R^+) \leq \omega(x^+ k_R^+) = 0, \quad \forall \ x^+ \in \Aa_L^+.
\end{equation}
Proposition~\ref{Pro:K2} then assures us that $a_R^+$ belongs to $\Kk_\omega \cap \Aa_R^+$.\qed

\vspace{0.2cm}

From above and Proposition~\ref{Pro:OrderQ}, we can conclude that $(\Bb_\omega,\Dd_1)$ is an order space with unit $1+\Kk_\omega$. As such, $\Bb_\omega$ can be endowed with an order seminorm $\llbracket \cdot \rrbracket$. Since the parent ordered space of $\Kk_\omega$ is a $C^\ast$-algebra, Proposition~\ref{Pro:KIdeal} has actually far more reaching consequences, as explained next.

\subsection{Background: Operator systems and matrix ordered $\ast$-vector spaces} We collect here a number of fundamental concepts and statements related to order structures on matrix amplifications.

\begin{definition}[\cite{PaulsenBook2002}~p.~9]\label{Def:OS1} A concrete operator system is a self-adjoint linear subspace of a unital $C^\ast$-algebra containing the unit.
\end{definition}

\begin{remark}{\rm As in the case of operator spaces, we call an operator system closed if the linear subspace in Definition~\ref{Def:OS1} is closed.
}$\Diamond$
\end{remark}

A concrete operator system inherits a full order structure from the embedding unital $C^\ast$-algebra. Indeed, if $\Ss \subseteq \Aa$ is an operator system, then $\Ss^+ = \Ss \cap \Aa^+$ supplies a positive cone. A matrix amplification of an operator system is again a linear subspace of a $C^\ast$-algebra, which is the matrix amplification of the embedding $C^\ast$-algebra. As such, the matrix amplifications of an operator system come equipped with order structures too. This tower of order structures puts a sharp distinction between the linear spaces that can or can not be embedded in $C^\ast$-algebras such that they contain the unit. The mentioned extra structures can be described abstractly and intrinsically.  

\begin{definition}[\cite{PaulsenBook2002}~p.~176]\label{Def:MatOrder} Given a $\ast$-vector space $\Vv$, one says that $\Vv$ is matrix-ordered provided that:
\begin{enumerate}[{\ \ \rm i)}]
\item For each $n$, we are given a cone $\Cc_n$ in $M_n(\Vv)_{\rm h}$;
\item $\Cc_n \cap (-\Cc_n) = \{0\}$ for all $n$;
\item For every $n,m\in \NM^\times$ and $A \in M_{n,m}(\CM)$, we have $A^\ast \Cc_n A \subseteq \Cc_m$.
\end{enumerate}
\end{definition}

\begin{definition} Let $(\Vv,\Vv^+)$ be a matrix-ordered $\ast$-vector space with order unit $e$. Then $e$ is called a matrix order unit provided $I_n \otimes e$ is an order unit for $M_n(\Vv)$, for each $n$. Furthermore, $e$ is called Archimedean matrix order unit if $I_n \otimes e$ is an Archimedean order unit for $M_n(\Vv)$, for each $n$.
\end{definition} 

\begin{definition}[\cite{PaulsenBook2002}~p.~176] Given two matrix-ordered $\ast$-vector spaces $\Vv$ and $\Vv'$ with cones $\Cc_n$ and $\Cc'_n$, one calls a linear map $\varphi : \Vv \rightarrow \Vv'$ completely positive provided that $[v_{ij}]\in \Cc_n$ implies that $[\varphi(v_{ij})] \in \Cc'_n$. One calls $\varphi$ a complete order isomorphism if $\varphi$ is completely positive and it has an inverse which is also completely positive.
\end{definition}

The following result, due to Choi and Effros \cite{ChiJFA1977}, supplies the abstract characterization of operator systems.

\begin{theorem}[\cite{PaulsenBook2002}~p.~177]\label{Th:CE} If $\Vv$ is a matrix-ordered $\ast$-vector space with an Archimedean matrix order unit $e$, then there exists a Hilbert space $H$, a concrete operator system $\Ss \in B(H)$, and a complete order isomorphism $\varphi : \Vv \rightarrow \Ss$ with $\varphi(e) = I_H$. Conversely, every concrete operator system $\Ss$ is a matrix-ordered $\ast$-vector space with Archimedean matrix order unit, when equipped with the matrix order inherited from the embedding $C^\ast$-algebra and with the Archimedean matrix order unit $e = 1$.
\end{theorem}

\begin{definition} A linear map between two abstract operator systems is called unital if it maps the order unit into the order unit. As in \cite{KavrukJFA2011}, we denote by $\Oo$ the category whose objects are operator systems and whose morphisms are unital completely positive (u.c.p.) maps. 
\end{definition}

The following statement supplies an effective criterion for a map to be completely positive. It will be used here very often.

\begin{proposition}[\cite{BlecherBook},~p.~18, \ \cite{PaulsenBook2002}~Prop.~2.11]\label{Pro:UnitalContr} Let $\Ss$ and $\Ss'$ be two operator systems and $\varphi : \Ss \rightarrow \Ss'$ be a linear unital map. Then $\varphi$ is completely positive if and only if $\varphi$ is completely contractive for the order norms.
\end{proposition}

\begin{remark}\label{Re:OsyVsOsp}{\rm Together with Remark~\ref{Re:CStarOrder}, Proposition~\ref{Pro:UnitalContr} implies that the unital complete order embedding of an abstract operator space $\Ss$ in a $B(H)$ supplies also isometric embeddings of $(M_n(\Ss), \llbracket \cdot \rrbracket_n)$ in $B(H^{(n)})$, for all $n \in \NM^\times$. This means $(\Ss, \llbracket \cdot \rrbracket)$ is an operator space and that its system of matrix norms $\|\cdot \|_n$ coincide with $\llbracket \cdot \rrbracket_n$. The conclusion is that $\llbracket \cdot \rrbracket_n$ satisfies Ruan's axioms.
}$\Diamond$
\end{remark}

Stinespring theorem \cite{StinespringAMS1955}, formulated below, supplies the structure of the completely positive maps when the domain is a $C^\ast$-algebra and the codomain is $B(H)$ of some Hilbert space.

\begin{theorem}[\cite{BlecherBook},~p.~18]\label{Th:Stinespring} Let $\Aa$ be a unital $C^\ast$-algebra. A linear map $\varphi : \Aa \rightarrow B(H)$ is completely positive if and only if there is a Hilbert space $K$, a unital $\ast$-homomorphism $\pi : \Aa \rightarrow B(K)$, and a bounded linear map $V :H \rightarrow K$ such that $\varphi(a)=V^\ast \pi(a) V$ for all $a \in \Aa$. This can be accomplished with $\|\varphi\|_{\rm cb} = \|V\|^2$. Also, this equals $\|\varphi\|$. If $\varphi$ is unital, then we may take $V$ to be an isometry; in this case we may view $H \subseteq K$ and have $\varphi(a)=P_H \pi(a) |_H$
\end{theorem}

Arveson extension theorem \cite{AversonAA1969}, formulated below, tells us among many other things that the above factorization functions also when the domain is an operator system.

\begin{theorem}[\cite{BlecherBook},~p.~18]\label{Th:Arveson} If $\Ss$ is an operator subsystem of a unital $C^\ast$-algebra $\Aa$, and if $\varphi: \Ss \rightarrow B(H)$ is completely positive, then there exists a completely positive map $\hat \varphi : \Aa \rightarrow B(H)$ extending $\varphi$.
\end{theorem}

\subsection{The reduced space is a matrix-ordered $\ast$-space} We have found in Proposition~\ref{Pro:KIdeal} that $\Kk_\omega$ is an order ideal. This, together with the fact that $\Kk_\omega$ is a subspace of a $C^\ast$-algebra enables us to formulate one of our main conclusions:

\begin{proposition} Let $\omega$ be any state on $\Aa^{\otimes \ZM}$ and recall the subsets $\Dd_n$ of $\Bb_\omega$, introduced in Eq.~\eqref{Eq:Dn}. Then the reduced data 
\begin{equation}
\big (\Bb_\omega, \{\Dd_n\}_{n \geq 1}, 1+\Kk_\omega \big )
\end{equation}
defines a matrix-ordered $\ast$-space with a matrix-order unit.
\end{proposition}

\begin{proof} We reproduce the discussion in \cite{KavrukAM2013}[p.~327] from where we learn that, for any linear self-adjoint subspace of an operator system that does not contain the unit, the projections of the positive cones are also cones for the matrix amplifications of the quotient space and, furthermore, they automatically satisfy the compatibility conditions iii) from Definition~\ref{Def:MatOrder}. Furthermore, if this linear subspace is an order ideal, condition ii) from Definition~\ref{Def:MatOrder} is also satisfied.
\end{proof}

Thus, the images of the positives cones of $\Aa_R$ and of its matrix amplifications through the quotient map $q : \Aa_R \to \Aa_R/\Kk_\omega$ and its matrix amplifications supply a matrix-ordered $\ast$-space structure with a unit on $\Bb_\omega$, and this happens for any state $\omega$ on $\Aa_\ZM$. Our next task is to complete the matrix-ordered structure to an Archimedean one, to confirm that $\Bb_\omega$ inherits from its parent $C^\ast$-algebra both, an operator space structure and an operator system structure. This will also help us elucidate the fate of $\omega$ under the reduction process.

\subsection{Archimedeanization of the reduced space} There are well understood Archimedeanization processes, which were developed in \cite{PaulsenIUMJ2009} for ordered vector spaces and in \cite{PaulsenPLMS2010} for matrix amplifications. They typically involve two stages, of which the first one quotients out the kernels of the matrix seminorms and the second one expands the positive cones. It was shown in \cite{KavrukAM2013} that this process greatly simplifies if the base space is a quotient by a kernel:

\begin{definition}[\cite{KavrukAM2013}, \ Def.~3.1] A subset $J$ of an operator system $\Ss$ is called a kernel if there exists a collection $\{\eta_\alpha\}_{\alpha \in \AM}$ of states on $\Ss$ such that $J = \bigcap_{\alpha \in \AM} {\rm ker} \, \eta_\alpha$.
\end{definition}

\begin{proposition}[\cite{KavrukAM2013}, \ Lemma~3.3] Let $J$ be a closed, non-unital order ideal of an operator system $\Ss$. Then the order seminorm on $\Ss/J$ is a norm if and only if $J$ is a kernel.
\end{proposition} 

Quite remarkably, for {\it any} state $\omega$ on $\Aa_\ZM$, the entanglement kernel is an order kernel ideal. We establish this fact in several steps.

\begin{proposition} Let $x \in \Aa_L^+$ such that $\omega(x)=0$. Then ${\rm Ker}\, \omega_x = \Aa_R$.
\end{proposition}

\begin{proof} By renormalizing $x\in \Aa_L^+$ by its norm in $\Aa_\ZM$, we can assume $\|x\|=1$. From the Cauchy-Schwarz inequality, we have:
\begin{equation}\label{Eq:CS}
|\omega(x a_R)|^2 \leq \omega(x^2) \, \omega(a_R^\ast a_R), \quad \forall \ a_R \in \Aa_R.
\end{equation}
Since $\|x\|\leq 1$, we have $(1-x)x \geq 0$ or $x^2 \leq x$, hence $0 \leq \omega(x^2) \leq \omega(x)=0$. Then \eqref{Eq:CS} assures us that $\omega_x(a_R)=0$ for all $a_R \in \Aa_R$.\end{proof}

\begin{corollary} The entanglement kernel can be equivalently defined as
\begin{equation}\label{Eq:Kappa3}
\Kk_\omega = \bigcap_{x \in \Aa_L^+}^{\omega(x) \neq 0} {\rm Ker} \, \omega_x.
\end{equation}
\end{corollary}

\begin{proof} Indeed,
\begin{align}
\bigcap_{x \in \Aa_L^+} {\rm Ker} \, \omega_x & = \Big ( \bigcap_{x \in \Aa_L^+}^{\omega(x) \neq 0} {\rm Ker} \, \omega_x \Big ) \bigcap \Big ( \bigcap_{x \in \Aa_L^+}^{\omega(x) = 0} {\rm Ker} \, \omega_x \Big ) \\ \nonumber
& = \Big ( \bigcap_{x \in \Aa_L^+}^{\omega(x) \neq 0} {\rm Ker} \, \omega_x \Big ) \bigcap \Aa_R
\end{align}
and the statement follows.\end{proof}

The value of the last statement rests in the observation that all positive functionals $\omega_x$ entering in the new definition \eqref{Eq:Kappa3} of $\Kk_\omega$ can be normalized by $\omega_x(1)$, hence, transformed into states. More precisely:

\begin{proposition} The entanglement kernel is a kernel. Explicitly, the entanglement kernel is the intersection of the kernels of a family of states:
\begin{equation}\label{Eq:Kappa4}
\Kk_\omega = \bigcap_{x \in \Aa_L^+}^{\omega(x) = 1} {\rm Ker} \, \omega_x.
\end{equation}
Note that $\omega_x$ is a state on $\Aa_R$ if $\omega(x) =1$.
\end{proposition} 

\begin{proof} We have
\begin{equation}
\bigcap_{x \in \Aa_L^+}^{\omega(x) \neq 0} {\rm Ker} \, \omega_x = \bigcap_{\alpha \in (0,\infty)} \Big (\, \bigcap_{x \in \Aa_L^+}^{\omega(x) =\alpha} {\rm Ker} \, \omega_x \Big ).
\end{equation}
Obviously, ${\rm Ker}\, \omega_x = {\rm Ker} \, \omega_{\alpha x}$ for all $\alpha \in (0,\infty)$, hence,
\begin{equation}
\bigcap_{x \in \Aa_L^+}^{\omega(x) =\alpha} {\rm Ker} \, \omega_x = \bigcap_{x \in \Aa_L^+}^{\omega(x) =1 } {\rm Ker} \, \omega_x, \quad \forall \ \alpha \in (0,\infty),
\end{equation}
because both, the left and right sides, sample the same subsets of $\Aa_R$. The statement then follows. \end{proof}

We arrive at the main result of the section: 

\begin{theorem}\label{Th:FinalB} The reduced space $\Bb_\omega = \Aa_R/\Kk_\omega$ inherits from $\Aa_R$ a canonical Archimedean matrix-order structure, supplied by the Archemedianization of $(\Bb_\omega, \{\Dd_n\}_{n \geq 1})$. Its positive cones are (see \cite{KavrukAM2013}[Prop.~3.4] and \cite{PaulsenPLMS2010}[Prop.~3.16])
\begin{equation}\label{Eq:CnS}
\begin{aligned}
\Cc_n = \{[a^{ij}_R+\Kk_\omega] \in M_n(\Bb_\omega) \ |& \ \forall \, \epsilon >0 \ \exists \, k_{ij} \in \Kk_\omega \ \mbox{such that} \ \\ 
& \quad \epsilon 1 \otimes I_n + [a^{ij}_R + k_{ij}] \in M_n(A_R)^+\big \},
\end{aligned}
\end{equation}
and its Archimedean matrix unit is $e=1+\Kk_\omega$. Furthermore, the quotient map $q:\Aa_R \rightarrow \Bb_\omega$ is u.c.p..
\end{theorem}

\begin{remark}{\rm One important point about the above is that $\Bb_\omega$ is not being quotiented because its order norm induced by $\Dd_1$ and unit $1+\Kk_\omega$ was already a norm. Note, however, that this does not imply that its unit was Archimedean.
}$\Diamond$
\end{remark}

The Archimedeanization of $\Bb_\omega$ we just described enjoys the universal property described in \cite{KavrukAM2013}[p.~329], which can be used to characterize the reduced map $\bar \omega$:

\begin{proposition}[\cite{KavrukAM2013}~Prop.~3.16] Assume that $\Bb_\omega$ is equipped with operator space structure as in Theorem~\ref{Th:FinalB}. Let $\Tt$ be an operator system and $\varphi : \Aa_R \rightarrow \Tt$ be a unital and completely positive map such that $\Kk_\omega \subseteq {\rm Ker}\, \varphi$. Then the map $\bar \varphi : \Bb_\omega \rightarrow \Tt$ given by $\bar \varphi(a_R + \Kk_\omega) = \varphi(a_R)$ is unital and completely positive. In particular, $\bar \omega$ is a unital and completely positive map from $\Bb_\omega$ to $\CM$, hence a state.
\end{proposition}

We, actually, can say much more:

\begin{proposition}\label{Pro:BarO1} $\bar \omega \big (\Dd_1 \setminus \{0\} \big ) \cap \{0\} = \emptyset$.
\end{proposition}

\begin{proof} Let $a_R \in \Aa_R^+$, which we can always normalize such that $0 \leq a_R \leq 1$. We will show that $\bar \omega(\hat a_R) =0$ implies $a_R \in \Kk_\omega$. Indeed, take any $x$ from $\Aa_L^+$. From Cauchy-Schwarz inequality, we have that $\omega(x a_R)^2 \leq \omega(x^2) \omega(a_R^2)$. Now, with the assumed normalization, $a_R(1-a_R) \geq 0$, or $a_R \geq a_R^2$. Then
\begin{equation}
\omega(x a_R)^2 \leq \omega(x^2) \omega(a_R^2) \leq \omega(x^2) \omega(a_R) =0,
\end{equation}
which proves that $\omega(x a_R) =0$ for all $x \in \Aa_L^+$, hence $a_R \in \Kk_\omega$.
\end{proof}

As in \cite{KavrukAM2013}, we denote by $\|\cdot \|_n^{\rm osy}$ the order norms on matrix space $M_n(\Bb_\omega)$, which obey Ruan's axioms (see Remark~\ref{Re:OsyVsOsp}). From the discussion in \cite{PaulsenPLMS2010}[p.~37], one learns that the $\Cc_n$'s defined above are just the closures of $\Dd_n$'s in the topology induced by these order norms. Furthermore, the quotient map $q:\Aa_R \rightarrow \Bb_\omega$ and its matrix amplifications are unital and contractive for $\|\cdot \|_n^{\rm osy}$. An important question we need to answer is if these norms are complete, which is equivalent to asking if the induced operator system structure is closed. Below is the answer:

\begin{proposition}\label{Pro:Complete} The norms $\|\cdot \|_n^{\rm osy}$ are complete if and only if they are equivalent with $\|\cdot \|_n^{\rm osp}$. If that is the case, then $\Dd_n$'s coincide with $\Cc_n$'s, hence, $\{\Bb_\omega,\{\Dd_n\}_{n \geq 1},e)$ is a \underline{closed} operator system. Furthermore, the reduced state $\bar \omega$ is faithful.
\end{proposition}

\begin{proof} As already pointed out in \cite{KavrukAM2013} (see also Corollary~\ref{Cor:42}), $\| \cdot \|_n^{\rm osy} \leq \|\cdot \|_n^{\rm osp}$ and we know that $\|\cdot \|_n^{\rm osp}$ are complete. If $\| \cdot \|_n^{\rm osy}$ is complete, then the identity map on $M_n(\Bb_\omega)$ is a continuous map from the Banach space $(M_n(\Bb_\omega), \| \cdot \|_n^{\rm osp})$ to the Banach space $(M_n(\Bb_\omega), \| \cdot \|_n^{\rm osy})$. Since the identity map is surjective, it is automatically an open map, which means its inverse, mapping $(M_n(\Bb_\omega), \| \cdot \|_n^{\rm osy})$ into $(M_n(\Bb_\omega), \| \cdot \|_n^{\rm osp})$, is continuous. This imply the existence of finite positive constant $k$ such that $\| \cdot \|_n^{\rm osp} \leq  k \| \cdot \|_n^{\rm osy}$. The second statement follows from Proposition~\ref{Pro:BarO1}.\end{proof}

In general, the operator space and operator system norms on quotient spaces are not equivalent (see \cite{KavrukAM2013}[Sec.~4]). As such, if one has a preference for closed operator systems, $\Bb_\omega$ needs to be completed. This can be done straightforwardly using a concrete representation. Indeed, complete or not, Choi-Effros Theorem~\ref{Th:CE} assures us that:

\begin{corollary}\label{Cor:Embedding} The quotient space with its Archimedean matrix-order structure $(\Bb_\omega,\{\Cc_n\}_{n \geq 1}, e)$ admits a concrete representation as an operator system inside $B(H)$ of some Hilber space $H$. The closure of that representation inside $B(H)$ supply the completion of $\Bb_\omega$.
\end{corollary}

Furthermore, from Stinespring \cite{StinespringAMS1955} and Arveson Theorems~\ref{Th:Stinespring} and \ref{Th:Arveson}, respectively, we now can spell out the structure of the reduced map $\bar \omega$:

\begin{corollary} If $\rho : \Bb_\omega \rightarrow B(H)$ is the concrete representation of $\Bb_\omega$, then there exists a vector $\zeta \in H$ such that
\begin{equation}
\bar \omega(b) = \big \langle \zeta, \rho(b) \zeta \big \rangle, \quad b \in \Bb_\omega.
\end{equation}
Furthermore, the map extends to a completely positive map over the entire $B(H)$ and, in particular, over the completion of $\Bb_\omega$.
\end{corollary}

\subsection{Sufficient conditions for completeness of order norms} We can spell out a simple and explicit condition that ensures that the matrix order norms are complete. This condition will appear again later in a totally different context.

\begin{proposition}[\cite{KavrukAM2013},~Prop.~4.1]\label{Pro:NormC} The operator space and the matrix order norms on $\Bb_\omega$ can be characterized as
\begin{align}
\|[a_R^{ij} + \Kk_\omega ]\|^{\rm osp}_n = \sup \big \{\| & [\phi(a_R^{ij})]\| \ |  \ \phi: \Aa_R \rightarrow B(H), \\ \nonumber
& \phi(\Kk_\omega) = \{0\},\ \phi \ \mbox{completely contractive} \big \}
\end{align}
and
\begin{align}
\| [a_R^{ij} + \Kk_\omega] \|^{\rm osy}_n = \sup \big \{\|[\phi(a_R^{ij})]\| \ | & \ \phi: \Aa_R \rightarrow B(H), \ \phi(\Kk_\omega) = \{0\}, \\ \nonumber
& \quad \phi \ \mbox{unital completely positive} \big \}
\end{align}
where, in each of the cases, $H$ runs through all Hilbert spaces.
\end{proposition} 

The following are direct consequences of the above:

\begin{corollary}[\cite{KavrukAM2013},~Cor.~4.2]\label{Cor:42} $\|\cdot \|^{\rm osy}_n \leq \|\cdot \|^{\rm osp}_n$.
\end{corollary}

\begin{corollary} The following inequality
\begin{equation}\label{Eq:Ineq2}
\| [a_R^{ij}+\Kk_\omega] \|^{\rm osy}_n \geq \sup_{x \in \Aa_L^+}^{\omega(x)=1} \|[\omega_x(a_R^{ij})]\|_{M_n(\CM)}
\end{equation}
 holds for any $[a_R^{ij}+\Kk_\omega] \in M_n(\Bb_\omega)$.
 \end{corollary}
 
 \begin{proof} This follows from Proposition~\ref{Pro:NormC} because $\omega_x$ is a state over $\Aa_R$ if $x \in \Aa_L^+$ is such that $\omega(x) =1$.\end{proof}
 
The above prompts us to define:

\vspace{0.2cm}

\begin{proposition}\label{Pro:Gamma} The map
\begin{equation}\label{Eq:Gamma}
\Gamma: \Bb_\omega \rightarrow [0,\infty), \quad \Gamma(a_R+\Kk_\omega) : = \sup_{x \in \Aa_L^+}^{\omega(x)=1} |\omega_x(a_R)|,
\end{equation}
and its generalizations
\begin{equation}\label{Eq:GammaN}
\Gamma_n: M_n(\Bb_\omega) \rightarrow [0,\infty), \quad \Gamma_n([a_R^{ij}+\Kk_\omega]) : = \sup_{x \in \Aa_L^+}^{\omega(x)=1} \|[\omega_x(a_R^{ij}) ] \|_{M_n(\CM)},
\end{equation}
for $n \in \NM^\times$, are well defined, continuous, sub-linear and homogeneous. Furthermore, 
\begin{equation}
\Gamma_n\big (M_n(\Bb_\omega) \setminus \{0\}\big) \subset (0,\infty).
\end{equation}
\end{proposition}

\begin{proof} The maps are well defined for, if $[b^{ij}_R]$ is another element from the class of $[a^{ij}_R]$, then $b^{ij}_R - a^{ij}_R \in {\rm Ker} \ \omega_x$ and $\omega_x(b^{ij}_R) = \omega_x(a^{ij}_R)$, for all $x \in \Aa_L^+$ with $\omega(x)=1$. Sub-linearity follows from the linearity of each $\omega_x$ and from the ``sub-linearity'' of $\| \cdot \|_{M_n(\CM)}$ and of the sup process. If $\alpha \in (0,\infty)$, then
\begin{equation}
\begin{aligned}
& \sup_{x \in \Aa_L^+}^{\omega(x)=1} \|[\omega_x(\alpha a^{ij}_R) ]\|_{M_n(\CM)} = \sup_{x \in \Aa_L^+}^{\omega(x)=1} \|  [\alpha \omega_x(a^{ij}_R) ]\|_{M_n(\CM)}  \\
& \qquad = \sup_{x \in \Aa_L^+}^{\omega(x)=1} \|\alpha[\omega_x( a^{ij}_R) ]\|_{M_n(\CM)} = \sup_{x \in \Aa_L^+}^{\omega(x)=1} \alpha \|[\omega_x( a^{ij}_R) ]\|_{M_n(\CM)}
\end{aligned}
\end{equation}
and homogeneity follows. Lastly, $[a^{ij}_R+\Kk_\omega] \in M_n(\Bb_\omega) \setminus \{0\}$ implies that $[a^{ij}_R] \notin M_n(\Kk_\omega)$, hence, there exists at least one $x \in \Aa_L^+$ with $\omega(x) = 1$ and a pair ${i,j}$ such that $\omega_x(a^{ij}_R) >0$. Then the last statement follows. \end{proof}

We are now ready to supply the sought condition:

\begin{proposition}\label{Pro:ArchCond} If the closure of the image of the unit ball of $(\Bb_\omega,\|\cdot \|^{\rm osp}_1)$ through the map $\Gamma$ does not contain the origin of the real axis, then $\|\cdot\|^{\rm osp}_n$ and $\|\cdot \|^{\rm osy}_n$ are all equivalent and $\Dd_n$ coincides with $\Cc_n$.
\end{proposition}

\begin{proof} With the stated assumption, there must exist a strictly positive constant $c$ such that $\Gamma(b)  \geq c$, for all $b \in \Bb_\omega$ with $\|b\|^{\rm osp}_1 = 1$. Since $\Gamma$ is homogeneous, this is same as $\Gamma(b) \geq c \, \|b\|^{\rm osp}_1$, for all $b \in \Bb_\omega$. This together with \eqref{Eq:Ineq2} and Corollary~\ref{Cor:42} give
\begin{equation}
\|b \|_1^{\rm osp} \geq  \|b \|_1^{\rm osy} \geq c \, \|b \|_1^{\rm osp}.
\end{equation}
This implies $\| \cdot \|_1^{\rm osy}$ is complete and, from Remark~\ref{Re:OsyVsOsp}, we known that $\| \cdot \|_n^{\rm osy}$ satisfy Ruan's axioms. Then Proposition~\ref{Pro:RuanVsNorms} assures us that all matrix norms $\| \cdot \|_n^{\rm osy}$ are complete, hence equivalent to $\| \cdot \|_n^{\rm osp}$. \end{proof}

\begin{remark}{\rm In the language introduced in \cite{KavrukAM2013}[p.~334], the entanglement kernel $\Kk_\omega$ becomes completely order proximinal under the conditions of Proposition~\ref{Pro:ArchCond}. It will be interesting to establish if these conditions are optimal for $\Kk_\omega$ to enjoy this property. One should be aware that these conditions are still not sufficient for the operator space and system structures to be isometric (see \cite{KavrukAM2013}[Prop.~4.10]).
}$\Diamond$
\end{remark} 

\begin{corollary}\label{Cor:FD} If $\Bb_\omega$ is finite dimensional, then $\| \cdot \|_n^{\rm osp}$ and $\| \cdot \|_n^{\rm osy}$ are equivalent.
\end{corollary}

\begin{proof} Any two norms on a finite linear space are equivalent.\end{proof}

\begin{remark}{\rm The reader should not be deceived by the simplicity of the above statement or of its proof. Indeed, they depend crucially on the fact that $\| \cdot \|_n^{\rm osy}$ are norms and not mere semi-norms and the latter follows from the deep inside from \cite{KavrukAM2013}, and the amazing fact that $\Kk_\omega$ is a kernel.}$\Diamond$
\end{remark}

\subsection{Concluding remarks and a look ahead} The work \cite{FannesCMP1992} exposed the matrix ordered $\ast$-space structure of $\Bb_\omega$ in Lemma~A1, under the assumptions that $\omega$ is translation invariant and $\Bb_\omega$ is finite dimensional. This, however, is not enough for $\| \cdot \|_n^{\rm osy}$ to be a norm and set Corollary~\ref{Cor:FD} in motion, or, equivalently, to embed $\Bb_\omega$ with the order structure specified in \cite{FannesCMP1992}[Lemma~A1]  in a $C^\ast$-algebra. As indicated in the discussion at page 451 of \cite{FannesCMP1992}, this is not at all a concerned if the minimality of $\Bb_\omega$ is not enforced. As we mentioned at the beginning of the section, the works of Kavruk, Paulsen, Todorov and Tomforde \cite{PaulsenIUMJ2009,PaulsenPLMS2010,KavrukJFA2011,KavrukAM2013} provide just the right tools and, in fact,  an entire framework to completely settle such fine points for any state over $\Aa_\ZM$ with $\Aa$ nuclear. What we have learned in this section is that, without truncation or completion, the quotient linear space $\Bb_\omega = \Aa_\ZM/\Kk_\omega$ inherits a canonical operator system structure from its parent $C^\ast$-algebra $\Aa_\ZM$, with the positive cones specified in Eq.~\eqref{Eq:CnS}. The latter can be different from the images of the positive cones of $\Aa_\ZM$ through the quotient maps and the induced order topology can also be different from the quotient topology.

Proposition~\ref{Pro:ArchCond} supplies sufficient conditions for the operator system to be closed. If $\Bb_\omega$ is finite dimensional, the unit ball of $(\Bb_\omega),\|\cdot \|^{\rm osp}_1)$ is compact and its image through $\Gamma$ is also compact, hence it cannot contain the origin of the real axis. Thus, the condition in Proposition~\ref{Pro:ArchCond} is automatically satisfied for finite dimensional reduced spaces and, as such, it is reasonable to claim that the situations singled out in Proposition~\ref{Pro:ArchCond} are direct generalizations of the cases studied in \cite{FannesCMP1992}. Certainly, they are very interesting cases to study in the future.

Looking forward, we need to make a choice for the reduced data. Our choice is the Archimedean matrix order $\ast$-vector space $\Bb_\omega$ without completion, together with its matrix order norms and the reduced state. Thus, the tuple
\begin{equation}\label{Eq:FinalReducedD}
(\Bb_\omega, \{\Cc_n\}_{n \geq 1}, \| \cdot \|_n^{\rm osy}, e = 1 +\Kk_\omega, \bar \omega)
\end{equation}
represents our derived reduced data.

\section{Factorization Process}
\label{Sec:Factor}

The main conclusions of the previous section apply to generic states $\omega$ on the algebra of physical observables $\Aa_\ZM$. In this section, however, we start by assuming that the state is shift-invariant, $\omega \circ S = \omega$. In these conditions, \cite{FannesCMP1992} defined a bi-linear map $\Aa_\ZM \times \Bb_\omega \to \Bb_\omega$, which proved to be of fundamental importance in the analysis of one dimensional spin systems, as we already emphasized in our introductory remarks. In this section, we investigate the properties of this bi-linear map and of its extensions to tensor products. As for the previous phase of our program, key to this phase is the identification of the natural framework to work in. If the program was to be advanced inside the category of operator spaces, then the Haagerup tensor product of operator spaces is the right fit because of its natural relation with multi-linear forms \cite{PisierBook}[Ch.~5]. However, we already made the decision to advance the program inside the category of operator systems and, as such, we will place our analysis in the framework developed by Kavruk, Paulsen, Todorov and Tomforde \cite{KavrukJFA2011}, which systematizes the tensor product structures for operator systems.

\subsection{Background: Tensor products of operator systems} Throughout, $\odot$ denotes the algebraic tensor product of linear spaces.

\begin{definition}[\cite{KavrukJFA2011},~p.~273]\label{Def:TP} Let $(\Ss,\{\Pp_N\}_{n \geq 1},e_1)$ and $(\Tt,\{\Qq_N\}_{n \geq 1},e_1)$ be operator systems. An operator system structure on $\Ss \odot \Tt$ is a family $\tau=\{\Cc_n\}_{n \geq 1}$, $\Cc_n \subseteq \Ss \odot \Tt$, satisfying:
\begin{enumerate}[\rm T1.]
\item $(\Ss \odot \Tt, \{\Cc_n\}_{n\geq 1},e_1 \otimes e_2)$ is an operator system denoted $\Ss \otimes_\tau \Tt$;
\item $\Pp_n \odot \Qq_m \subset \Cc_{nm}$, for all $n,m \in \NM^\times$;
\item If $\phi: \Ss \to M_n(\CM)$ and $\psi : \Tt \to M_m(\CM)$ are u.c.p. maps, then $\phi \odot \psi : \Ss \otimes_\tau \Tt \to M_{nm}(\CM)$ is a u.c.p. map.
\end{enumerate}
\end{definition}

\begin{definition} $\otimes_\tau$ is called functorial if it can be extended to a functor from $\Oo \times \Oo$ to $\Oo$. 
\end{definition}

\begin{remark}\label{Re:TON}{\rm Given two operator system structures on $\Ss \odot \Tt$, one says that $\tau_1$ is greater than $\tau_2$ if the identity map on $\Ss \odot \Tt$ is a completely positive map from $\Ss \otimes_{\tau_1} \Tt$ to $\Ss \otimes_{\tau_2} \Tt$. This is equivalent to $M_n(\Ss \otimes_{\tau_1} \Tt)^+ \subseteq M_n(\Ss \otimes_{\tau_2} \Tt)^+$ for all $n \in \NM^\times$. In such case, the order norms enter the relation $\llbracket \cdot \rrbracket_{\tau_2} \leq \llbracket \cdot \rrbracket_{\tau_1}$.
}$\Diamond$
\end{remark}

Ref.~\cite{KavrukJFA2011} identified one side of the spectrum of operator tensor structures to be:

\begin{definition}[\cite{KavrukJFA2011},~p.~276] For $\Ss$ an operator system, let
\begin{equation}
S_n(\Ss) := \{\phi : \Ss \to M_n(\CM), \ \phi \ \mbox{\rm unital completely positive map}\}.
\end{equation}
The minimal tensor product $\Ss \otimes_{\rm min} \Tt$ of two operator systems $\Ss$ and $\Tt$ is defined by the system of positive cones
\begin{equation}
\begin{aligned}
\Cc_n^{\rm min}(\Ss,\Tt)  : =  \{[p_{ij}] \in M_n(\Ss \odot \Tt), & \ [(\phi \odot \psi) (p_{ij})] \in M_{nkm}(\CM)^+, \\
& \quad  \phi \in S_k(\Ss), \ \psi \in S_m(\Tt), \ k,m\in \NM^\times\}.
\end{aligned}
\end{equation}
\end{definition}

\begin{theorem}[\cite{KavrukJFA2011},~Th.~4.4] Let $\Ss$ and $\Tt$ be operator systems and let $i_\Ss: \Ss \to B(H)$ and $i_\Tt : \Tt \to B(K)$ be embeddings that are unital complete order isomorphisms onto their ranges. Then ${\rm min}$ is the operator system structure on $\Ss \odot \Tt$ arising from the embedding $i_\Ss \odot \i_\Tt : \Ss \odot \Tt \to B(H \otimes K)$.
\end{theorem}

\begin{theorem}[\cite{KavrukJFA2011},~Th.~4.6] The mapping ${\rm min} : \Oo \times \Oo \to \Oo$ sending $(\Ss,\Tt)$ to $\Ss \otimes_{\rm min} \Tt$ is an injective, associative, symmetric, functorial operator system tensor product. Moreover, if $\tau$ is any other operator system structure on $\Ss \odot \Tt$, the $\tau$ is larger then ${\rm min}$.
\end{theorem}

On the other side of the spectrum sits:

\begin{definition}[\cite{KavrukJFA2011},~p.~276] The maximal tensor product $\Ss \otimes_{\rm max} \Tt$ of two operator systems $\Ss$ and $\Tt$ is defined by the Archimedeanization of the following system of positive cones:
\begin{equation}
\begin{aligned}
\Dd_n^{\rm max}(\Ss,\Tt) : =  \{\gamma ([s_{ij}] & \odot [t_{ij}]) \gamma^\ast,  \ [s_{ij}] \in M_k(\Ss)^+, \\
&  \ [t_{ij}] \in M_m(\Tt), \ \gamma \in M_{n,km}(\CM), \ k,m \in \NM^\times \}.
\end{aligned}
\end{equation}
\end{definition}

\begin{theorem}[\cite{KavrukJFA2011}, ~Th.~5.5] The mapping ${\rm max} : \Oo \times \Oo \to \Oo$ sending $(\Ss,\Tt)$ to $\Ss \otimes_{\rm max} \Tt$ is a symmetric, associative, functorial operator system product. Moreover, if $\tau$ is any other operator system structure on $\Ss \odot \Tt$, then ${\rm max}$ is larger than $\tau$.
\end{theorem}

\begin{remark}{\rm From \cite{KavrukJFA2011}[~Lemma~5.1], we learn that the matrix order induced by the positive cones $\Dd_n^{\rm max}$ is larger than any other matrix order satisfying property T2, in particular, it is larger than $\Ss \otimes_{\min} \Tt$. As such, Remark~\ref{Re:TON} assures us that the order seminorm induced by $\Dd_n^{\rm max}$ majorizes $\llbracket \cdot \rrbracket_{\rm min}$, hence, it is actually a norm. In this case, the Archimedeanization process consists of extending  $\Dd_n^{\rm max}$ to
\begin{equation}
\Cc_n^{\rm max}(\Ss,\Tt) := \{ [p_{ij}] \in M_n(\Ss \odot \Tt), \ r e_n + [p_{ij}] \in \Dd_n^{\rm max} \ \forall \ r >0\},
\end{equation}
and the first stage of the Archimedeanization process is not necessary.
}$\Diamond$
\end{remark}

\begin{theorem}[\cite{KavrukJFA2011},~Corollary 6.8]\label{Th:NuclearTP} Let $\Qq$ be a unital $C^\ast$-algebra. Then $\Qq$ is a nuclear $C^\ast$-algebra if and only if $\Qq \otimes_{\rm min} \Ss = \Qq \otimes_{\rm max} \Ss$ for every operator system $\Ss$.
\end{theorem}

The above statement assures us that $\Aa \odot \Bb_\omega$ carries only one operator structure structure if $\Aa$ is nuclear, which we denote simply by $\Aa \otimes \Bb_\omega$.

\subsection{Generating bi-linear map} The setting here is the same as in section~\ref{Sec:PhysAlg} but with the major difference that $\omega$ is assumed shift-invariant. In this section, the class $a_R + \Kk_\omega$ in $\Bb_\omega$ is denoted by $\hat a_R$ and $\Bb_\omega$ is considered equipped with its canonical operator system structure summarized in Eq.~\eqref{Eq:FinalReducedD}. We will continue to denote the quotient map from $\Aa_R$ to $\Bb_\omega$ by $q$.

Consider the $C^\ast$-algebra embeddings 
\begin{equation}\label{Eq:BetaMaps}
\Aa^{\otimes (p+1)} \rightarrowtail \Aa_\ZM, \quad \beta_{(n,n+p)} := \mathfrak i_{(n,n+p)} \circ (\otimes_{j=n}^{n+p} \, \alpha_j),
\end{equation}
for $p \geq 0$ and $n \in \ZM$.

\begin{proposition}\label{Pro:KO} Let $k \in \Kk_\omega \subseteq \Aa_R \subset \Aa_\ZM$. Then 
\begin{equation}
S^{\circ (p+1)} \big( \beta_{(-p,0)}(\alpha) k\big) \in \Kk_\omega,
\end{equation} 
for any $\alpha \in \Aa^{\otimes (p+1)}$ and $p \in \NM$.
\end{proposition}

\begin{proof} Let $x \in \Aa_L$. Given the shift invariance of $\omega$, we have 
\begin{equation}
\omega_x \big (S^{\circ (p+1)} ( \beta_{(-p,0)}(\alpha) k)\big ) = \omega\big(S^{\circ (-p-1)}(x) \beta_{(-p,0)}(\alpha) k \big).
\end{equation}
Therefore, if we denote by $x'$ the element $S^{\circ (-p-1)}(x) \beta_{(-p,0)}(\alpha)$ and observe that $x' \in \Aa_L$, then
\begin{equation}
\omega_x \big (S^{\circ (p+1)} ( \beta_{(-p,0)}(\alpha) k)\big ) = \omega_{x'}( k ) = 0
\end{equation}
and the statement follows.
\end{proof}

\begin{corollary} We have well-defined bilinear maps $\Aa^{\otimes (p+1)}  \times \Bb_\omega \to \Bb_\omega$,
\begin{equation}\label{Eq:EDef}
\EM_\omega^{(p+1)}(\alpha,\hat a_R) := (q\circ S^{\circ (p+1)})\big(\beta_{(-p,0)}(\alpha)a_R\big),
\end{equation}
for all $p \in \NM$.
\end{corollary}

We denote the canonical extension of $\EM_\omega^{(p+1)}$ as a linear map on the algebraic tensor product $\Aa^{\otimes (p+1)}  \odot \Bb_\omega$ by the same symbol and, for $p=0$, we simplify the notation to $\EM_\omega$. Since all operator system structures on $\Aa^{\otimes (p+1)}  \odot \Bb_\omega$ coincide, we have the liberty to choose between the various characterizations of this structure. For our purposes, we found that the max-structure is most efficient: 

\begin{proposition} $\EM_\omega^{(p+1)} : \Aa^{\otimes (p+1)}  \otimes_{\rm max} \Bb_\omega \to \Bb_\omega$ is a completely positive map.
\end{proposition}

\begin{proof} According to \cite{KavrukJFA2011}[Lemma~2.5], it is enough to check if 
\begin{equation}
\EM_{\omega,n}^{(p+1)}\big(\Dd_n^{\rm max}(\Aa^{\otimes (p+1)},\Bb_\omega)\big) \subseteq M_n(\Bb_\omega)^+.
\end{equation} For $[\alpha_{ij}] \in M_k(\Aa^{\otimes (p+1)})^+$ and $[\hat a_R^{ij}] \in M_m(\Bb_\omega)^+$, we have
\begin{equation}
\EM_{\omega,n}^{(p+1)} \big (\gamma ([\alpha_{ij}] \odot [\hat a_R^{ij}]) \gamma^\ast \big ) = \gamma \,  \EM_{\omega,km}^{(p+1)}([\alpha_{ij}] \odot [\hat a_R^{ij}]) \gamma^\ast,
\end{equation}
for any $\gamma \in M_{n,km}(\CM)$, and, by definition,
\begin{equation}
\EM_{\omega,km}^{(p+1)}([\alpha_{ij}] \odot [\hat a_R^{ij}]) = (q_{km} \circ S_{km}^{\circ (p+1)}) \big([\beta_{ij}] \otimes [a_R^{ij}]\big),
\end{equation}
where $[\beta_{ij}]=[\beta_{(-p,0)}(\alpha_{ij})]\in M_k(\Aa_L)^+$. Note that, on the right, we passed to the tensor product of $C^\ast$-algebras. According to Theorem~\ref{Th:FinalB}, our task is to show that, for any $\epsilon >0$, we can find $[K_{ij}] \in M_{km}(\Kk_\omega)$ such that
\begin{equation}
\epsilon e_{km}+ S_{km}^{\circ (p+1)} \big([\beta_{ij}] \otimes [a_R^{ij}]\big) +[K_{ij}]\in M_{km}(\Aa_R)^+.
\end{equation}
Since $[\hat a_R^{ij}]$ belongs to $M_m(\Bb_\omega)^+$, we know that, for any $\eta >0$, there exists $[k_{ij}(\eta)] \in M_m(\Kk_\omega)$ such that
\begin{equation}
\eta e_m + [a_R^{ij} + k_{ij}(\eta)] \in M_m(\Aa_R)^+.
\end{equation}
Then, if $\xi = \llbracket \alpha_{ij} \rrbracket_k = \llbracket \beta_{ij} \rrbracket_k \neq 0$, we choose
\begin{equation}
[K_{ij}] = S_{km}^{\circ (p+1)} \big([\beta_{ij}] \otimes [k_{ij}(\epsilon/\xi]\big),
\end{equation}
which is known to belong to $M_{km}(\Kk_\omega)$ by Proposition~\ref{Pro:KO}. We have
\begin{equation}
\begin{aligned}
& \epsilon e_{km} + S_{km}^{\circ (p+1)} \big([\beta_{ij}] \otimes [a_R^{ij}]\big) +[K_{ij}] \\
&\qquad = \epsilon e_{km} - \tfrac{\epsilon}{\xi} S_{km}^{\circ (p+1)} \big([\beta_{ij}] \otimes e_m\big) \\
 & \qquad \qquad + S_{km}^{\circ (p+1)} \big([\beta_{ij}] \otimes ( \tfrac{\epsilon}{\xi} e_m + [a_R^{ij}+k_{ij}(\epsilon/\xi)]\big),
\end{aligned}
\end{equation}
and the elements seen in the last two lines belong to the positive cone $M_{km}(\Aa_R)^+$. The case $\xi=0$ is evident.
\end{proof}

\begin{example}\label{Ex:ProdState2}{\rm For the product state from Example~\ref{Ex:ProdState1}, we found that $\Bb_\omega = \CM \cdot e$ and we have $\EM(a \otimes e) = \omega_0(a) \cdot e$.}$\Diamond$
\end{example}

\begin{proposition}\label{Pro:IsIter} The maps satisfy the recursive relations
\begin{equation}\label{Eq:Recursive}
\EM_\omega^{(p+1)} = \EM_\omega \circ ({\rm id} \otimes \EM_\omega^{(p)}),  \quad p > 1,
\end{equation}
and, more general,
\begin{equation}\label{Eq:Markov}
\EM_\omega^{(q+p)} = \EM_\omega^{(q)} \circ ({\rm id} \otimes \EM_\omega^{(p)}), \quad p ,q > 1.
\end{equation}
\end{proposition}

\begin{proof} The statements follow from the identity
\begin{equation}
S^{\circ (q+p)} \big(\beta_{(-q-p+1,0)}(\alpha \otimes \alpha') a_R\big) = S^{\circ q}\Big(\beta_{(-q+1,0)}(\alpha)  S^{\circ p}\big(\beta_{(-p+1,0)}(\alpha') a_R\big)\Big ),
\end{equation}
valid for any $\alpha \in \Aa^{\otimes q}$, $\alpha' \in \Aa^{\otimes p}$ and $a_R \in \Aa_R$.
\end{proof}

\begin{proposition}\label{Pro:IndL} Let us consider the unital complete order embeddings\footnote{The operator system structures on the tensor products are injective.} $J_p:\Aa^{\otimes p} \rightarrowtail \Aa^{\otimes p} \otimes \Bb_\omega$ and
$\mathfrak j_{p+k,p} : \Aa^{\otimes p} \rightarrowtail \Aa^{\otimes p} \otimes \Aa^{\otimes k} \simeq \Aa^{\otimes (p+k)}$.
Then
\begin{equation}\label{Eq:DirecL}
\EM_\omega^{(p+k)} \circ J_{p+k} \circ \mathfrak j_{p+k,p} = \EM_\omega^{(p)} \circ J_p
\end{equation}
on $\Aa^{\otimes p}$. As such, the tower of maps 
\begin{equation}
\EM_\omega^{(p)} \circ J_{p} : \Aa^{\otimes p} \to \Bb_\omega, \quad (p \geq 1),
\end{equation} 
has a direct limit $\EM_\omega^\infty \circ J_\infty: \Aa^{\otimes \NM^\times} \simeq \Aa_R \to \Bb_\omega$, which coincides with the quotient map $q: \Aa_R \to \Bb_\omega$.
\end{proposition}

\begin{proof} Eq.~\eqref{Eq:DirecL} is a consequence of the identity
\begin{equation}
S^{\circ p}\big(\beta_{(-p+1,0)}(\alpha) 1_R\big) = S^{\circ (p+k)}\big(\beta_{(-p+k+1,0)}(\alpha \otimes 1^{\otimes k}) 1_R\big).
\end{equation}
Eq.~\eqref{Eq:DirecL} then assures us that the tower of contractive unital maps respect the structure maps of the direct limit $\Aa^{\otimes \NM^\times}$, hence $\EM_\infty \circ J_\infty$ is well defined and shares the same attributes. Lastly, 
\begin{equation}
(\EM_\omega^{(p)}\circ J_p)(\alpha) = q \big( \beta_{(1,1+p)}(\alpha)\big ), \quad \forall \ \alpha \in \Aa^{\otimes p},
\end{equation} 
which proves the last statement.\end{proof}

\subsection{Ergodic states}
\label{Sub:Asy}

Shift-invariant states form a convex subset of the state space of $\Aa_\ZM$. The extremals among shift-invariant states, {\it i.e.} those states that can not be decomposed among the shift-invariant states, 
are called ergodic states \cite[Sec.~4.3]{BratteliBook1}. In this subsection, we investigate the synergy between the ergodic states and a dynamical system on $\Bb_\omega$, canonically induced by the shift-invariant state $\omega$. In particular, we supply sufficient conditions for the state $\omega$ to be ergodic.

Let $\pi_\omega : \Aa_\ZM \to B(H_\omega)$ be the GNS-representation induced by the shift-invariant state $\omega$. Then it is known that shift map $S$ is implemented by a conjugation with a unitary operator $U_\omega$, $\pi_\omega\big (S(\alpha)\big ) = U_\omega \pi_\omega(\alpha) U_\omega^\ast$. The following is a useful though somewhat abstract characterization of the ergodic states:

\begin{proposition}[\cite{BratteliBook1},~Th.~4.3.17]\label{Pro:Ergodic1} A shift-invariant state $\omega$ is ergodic if and only if the $C^\ast$-algebra generated by $\pi_\omega(\Aa_\ZM) \cup U_\omega$ inside $B (H_\omega)$ is irreducible on $H_\omega$.
\end{proposition}

A more practical criterion to identify ergodic states relies on asymptotic tests:

\begin{definition}\label{Def:Cluster1} We say that the shift-invariant state $\omega$ on $\Aa_\ZM$ has the asymptotic clustering property if the sequence
\begin{equation}\label{Eq:Cluster1}
\sup\Big \{\big | \omega \big (a_L \cdot S^{\circ r}(a_R) \big ) - \omega(a_L) \omega(a_R) \big |, \ a_{L} \in \Aa_{L}\subset \Aa_\ZM, \ \|a_{L}\| =1 \Big \}
\end{equation}
converges to zero as $r \to \infty$, for any $a_R \in \Aa_R\subset \Aa_\ZM$.
\end{definition}

\begin{proposition}[\cite{BratteliBook1},~Th.~4.3.22] If a shift-invariant state displays the asymptotic cluster property then the state is ergodic.
\end{proposition}

\begin{remark}{\rm The formulation of clustering property in Definition~\ref{Def:Cluster1} seems stronger than the standard formulation (see \cite[Sec.~4.3.2]{BratteliBook1}), in that it requires  a uniform convergence w.r.t. the $a_L$ entry. Note, however, that $a_L$ and $a_R$ are constrained on opposite half-sides of the chain in \ref{Eq:Cluster1}, in which case the formulation becomes equivalent with the standard one (see also Proposition~\ref{Pro:SMix}). The reason for our preference towards the formulation in Definition~\ref{Def:Cluster1} will become apparent in Proposition~\ref{Pro:FixPoint}.
}$\Diamond$
\end{remark}

From Proposition~\ref{Pro:Ergodic1}, we see that factor states are ergodic. In such cases, we can demonstrate that our formulation of clustering property can be indeed derived:

\begin{proposition}\label{Pro:SMix} If $\pi_\omega(\Aa_\ZM)''$ is a factor, then $\omega$ satisfies the cluster property as formulated in Definition~\ref{Def:Cluster1}.
\end{proposition} 

\begin{proof} We will follow closely the example~4.3.24 in \cite{BratteliBook1}, with a few improvements. The proof rests on the observation that $\Aa_\ZM$ displays the asymptotic abelienness 
\begin{equation}
\lim_{r \to \infty} \| [S^{\circ (-r)}(\alpha), \beta] \|=0, \quad \forall \ \alpha,\beta \in \Aa_\ZM,
\end{equation}
in a uniform fashion, provided $\alpha = a_L$ is chosen from $\Aa_L$. Indeed, for any $\epsilon>0$, there exist $N$ and $a_{N} \in \Aa_{[-N,N]}$ such that $\| \beta - a_{N}\| \leq \epsilon$. Then
\begin{equation}
\begin{aligned}
\| [S^{\circ (-r)}(a_L), \beta] \| = \| [S^{\circ (-r)}(a_L),\beta- a_{N}] \| \leq  2 \epsilon \|a_L\|,
\end{aligned}
\end{equation}
for all $r > N+1$. As a consequence,
\begin{equation}\label{Eq:Abel1}
\lim_{r \to \infty} \sup \big \{ \| [S^{\circ (-r)}(a_L), \beta]\|, \ a_L \in \Aa_L, \ \|a_L\|=1\big \} = 0.
\end{equation}
From here on, we can repeat the arguments from \cite[Example~4.3.24]{BratteliBook1}, which we do for completeness. Let $\Omega_\omega$ be the cyclic vector of $\pi_\omega$. Since $\Omega_\omega$ and $\eta=(\pi_\omega(a_R) - \omega(a_R) I) \Omega_\omega$ are orthogonal in $H_\omega$, there exists a self-adjoint operator $T$ on $H_\omega$ such that $T\Omega_\omega = \Omega_\omega$ and $T\eta = 0$. Taking $a_R$ self-adjoint and with $C=(\pi_\omega(a_R) - \omega(a_R) I)T$, we have $C \Omega_\omega = \eta$ and $C^\ast \Omega_\omega =0$, as well as
\begin{equation}
\begin{aligned}
\omega\big (a_L S^{\circ r}(a_R)\big ) - \omega(a_L) \omega(a_R)&  = \omega\big (S^{\circ (-r)}(a_L) (a_R-\omega(a_R) 1) ) \\
& = \langle \Omega_\omega, [\pi_\omega(S^{\circ (-r)}(a_L)),C]\Omega_\omega\rangle.
\end{aligned}
\end{equation}
Next, one observes that, since $\pi_\omega(\Aa^{\otimes \ZM})''$ is a factor, the algebra generated by $\pi_\omega(\Aa_\ZM) \cup \pi_\omega(\Aa_\ZM)'$ is irreducible on $H_\omega$, hence Kadison's transitivity theorem applies and $C$ can be chosen from this algebra. In particular, for each $\epsilon >0$, there exist $\beta_i \in \Aa_\ZM$ and $B_i \in \pi_\omega(\Aa_\ZM)'$, with $i$ in a finite set and such that
\begin{equation}
\| C - \sum\nolimits_i \pi_\omega(\beta_i) B_i \| \leq \epsilon.
\end{equation}
Then
\begin{equation}
\begin{aligned}
\langle \Omega_\omega, [\pi_\omega(S^{\circ (-r)}(a_L)),C]\Omega_\omega\rangle 
\leq 2 \epsilon \|a_L\|+ \sum_i \|B_i\| \| [S^{\circ (-r)}(a_L),\beta_i]\|
\end{aligned}
\end{equation}
and the statement follows from Eq.~\eqref{Eq:Abel1}. The case when $a_R$ is not self-adjoint is obvious.
\end{proof} 

We now investigate how ergodicity is related to the characteristics of a dynamical system on $\Bb_\omega$ canonically induced by the state $\omega$. Indeed, the shift map descends on $\Bb_\omega$:

\begin{definition} We call the shift map on $\Bb_\omega$ the map 
\begin{equation}
\bar S:\Bb_\omega \rightarrow \Bb_\omega, \quad \bar S = \EM_\omega \circ \bar L,
\end{equation} 
with  $\bar L$ being the unital complete order embedding
\begin{equation}
\bar L : \Bb_\omega \rightarrow \Aa \otimes \Bb_\omega, \quad \bar L(b) = 1 \otimes b.
\end{equation}
\end{definition}

\begin{proposition} As a composition of u.c.p. maps, $\bar S$ is u.c.p. and, furthermore, it satisfies the following relation $\bar S \circ q = q \circ S_R$. As a consequence, the reduced state is also shift invariant, $\bar \omega \circ \bar S = \bar \omega$.
\end{proposition}

\begin{proposition}\label{Pro:Cluster2} The shift-invariant state $\omega$ displays the asymptotic clustering property if and only if the reduced state $\bar \omega$ displays a similar reduced clustering property,
\begin{equation}\label{Eq:Cluster2}
\lim_{r \rightarrow \infty}\sup \Big \{ \Big | \bar \omega\Big(\EM_\omega^{(p)}\big (\alpha \otimes \bar S^{\circ r}(b) \big ) \Big ) - \bar \omega\Big (\EM_\omega^{(p)} \big (\alpha \otimes e\big )\Big ) \bar \omega(b)\Big | \Big \}=0,
\end{equation}
where the supremum is over $p \geq 1$ and $\alpha \in \Aa^{\otimes p}$ with $\|\alpha\|=1$.
\end{proposition}

\begin{proof} If $\beta = \beta_{(-p+1,0)}(\alpha) \in \Aa_L \subset \Aa_\ZM$ and $b=\hat a_R$ for some $a_R \in \Aa_R \subset \Aa_\ZM$, we have
\begin{equation}
\EM_\omega^{(p)} \big (\alpha \otimes \bar S^{\circ r}(b) \big )= q\big ( S^{\circ p} \big (\beta S^{\circ r}(a_R) \big ) \big ).
\end{equation}
Therefore, taking into account the shift invariance of $\omega$,
\begin{equation}
\bar \omega \big(\EM_\omega^{(p)} \big (\alpha \otimes \bar S^{\circ r}(b) \big ) \big ) = \omega\big (\beta \otimes S^{\circ r}(a_R) \big ).
\end{equation}
Then
\begin{align}
& \bar \omega\Big(\EM_\omega^{(p)} \big (\alpha \otimes (\bar S^{\circ r}(b) -\bar \omega(b) e )\big ) \Big )   = \omega\big (\beta \otimes S^{\circ r}(a_R) \big ) - \omega\big (\beta\big ) \omega\big ( a_R \big )
\end{align}
and, from this identity, the statement follows in both directions because $\bigcup_p \beta_{[-p,0]}(\Aa^{\otimes p})$ is dense in $\Aa_L$. \end{proof}

\begin{remark}\label{Re:CorrF}{\rm The asymptotic clustering property ensures the decay of correlation functions. For example, a 2-site correlation function refers a quantity of the type
\begin{equation}
C(a,a'; r) = \omega(\cdots \otimes 1 \otimes a \otimes 1^{\otimes r} \otimes a' \otimes 1 \otimes \cdots), \quad a,a' \in \Aa.
\end{equation}
This correlation function can also computed as 
\begin{equation}
C(a,a';r) = \bar \omega\circ \EM_\omega^{(r+2)}(a \otimes 1^{\otimes r} \otimes a' \otimes e) = (\bar \omega\circ \EM_\omega)\big (a \otimes \bar S^{\circ r}(\hat a') \big),
\end{equation}
and the clustering property assures us that $C(a,a';r) -\omega(a)\omega(a')$ converges to zero as the ``distance'' $r$ goes to infinity.
}$\Diamond$
\end{remark}

The next statement identifies specific conditions in which we can establish a direct relation between the dynamical system $(\Bb_\omega,\bar S)$ and the asymptotic clustering property of the state. For this, we recall the functional $\Gamma : \Bb_\omega \to [0,\infty)$ defined in Propostion~\ref{Pro:Gamma}. Then: 

\begin{proposition}\label{Pro:FixPoint} Assume  that
\begin{equation}\label{Eq:BarG}
\Gamma(b) \geq c \|b\|_1^{\rm osp}, \quad \ \forall \ b \in \Bb_\omega,
\end{equation} 
for some strictly positive constant $c$. Then the following are equivalent:
\begin{enumerate}[\ \rm 1.]
\item $\omega$ displays the asymptotic clustering property.
\item The linear sub-space generated by $e$ is the only attractor of the map $\bar S$,
\begin{equation}\label{Eq:FixP}
 \lim_{r \rightarrow \infty} \bar S^{\circ r} (b)= \bar \omega(b) \, e, \quad \forall \ b \in \Bb_\omega,
 \end{equation}
 where the limit is in the order topology of $\Bb_\omega$.
\end{enumerate}
\end{proposition}

\proof ``1 $\Rightarrow$ 2'' From the asymptotic clustering property~\eqref{Eq:Cluster1}, we have
\begin{equation}\label{Eq:X13}
\lim_{r \rightarrow \infty} \sup\Big \{ \Big | \omega\Big ( x \cdot S^{\circ r} \big (a_R - \omega(a_R) \, 1_R\big )  \Big ) \Big |, \ x \in \Aa_L \subset \Aa_\ZM, \ \|x\|=1 \Big \}= 0,
\end{equation}
for all $a_R \in \Aa_R \subset \Aa_\ZM$, which translates to
\begin{equation}
\Gamma \big ( \bar S^{\circ r}   ( \hat a_R - \omega(a_R) \, e  ) \big )  \rightarrow 0 \ \ {\rm as} \ \ r \rightarrow \infty.
\end{equation}
Since the unit is invariant for $\bar S$ and $\| \cdot \|_1^{\rm osp} \geq \| \cdot \|_1^{\rm osy}$, condition \eqref{Eq:BarG} implies
\begin{equation}
\big \| \bar S^{\circ r} (\hat a_R) - \bar \omega(\hat a_R) \, e \big\|_1^{\rm osy} \rightarrow 0 \ \ {\rm as} \ \ r \rightarrow \infty,
\end{equation}
which proves the first claim. 

``2 $\Rightarrow$ 1'' This is a direct consequence of Proposition~\ref{Pro:Cluster2}.\qed

\begin{remark}{\rm The above statement can be regarded as a direct generalization of point (3) of Proposition~3.1 in \cite{FannesCMP1992}.
}$\Diamond$
\end{remark}

\begin{corollary} In the conditions of Proposition~\ref{Pro:FixPoint}, there is one and only one shift-invariant state on $\Bb_\omega$, which can be detected via~\eqref{Eq:FixP}.
\end{corollary}

\begin{remark}{\rm Note that \eqref{Eq:BarG} is exactly the same condition that ensures that the operator system norm on $\Bb_\omega$ is complete (see Proposition~\ref{Pro:ArchCond}).
}$\Diamond$
\end{remark}

\begin{proposition}\label{Pro:EDense}  Assume that the statement in Eq.~\eqref{Eq:FixP} holds. Then the map $\EM_\omega$ is full, in the sense that
\begin{equation}\label{Eq:Full}
\overline{\bigcup\nolimits_{p \geq 1} \bigcup\nolimits_{x \in \Aa^{\otimes p}} \EM_\omega^{(p)}(x \otimes b)} = \Bb_\omega, \quad \forall \ b \in \Bb_\omega.
\end{equation}
\end{proposition}

\begin{proof} The statement is true for $b = e$ because $\EM_\omega^\infty \circ J_\infty = q$, as we have seen in Proposition~\ref{Pro:IndL}. Now, fix a generic element $b$ and let $\alpha_p \in \Aa^{\otimes p}$ be a uniformly bounded sequence such that $\EM_\omega^{(p)}(\alpha_p \otimes e)$ converges in order topology to another element $b' \in \Bb_\omega$ as $p \to \infty$. Then
\begin{equation}\label{Eq:X101}
\EM_\omega^{(2p)} \big ( (\alpha_p \otimes 1^{\otimes p}) \otimes b \big )  = \EM_\omega^{(p)}\big (\alpha_p \otimes \bar S^{\circ p}(b)\big ) .
\end{equation}
Rewriting the right side as
\begin{equation}\label{Eq:X102}
 \EM_\omega^{(p)}\big (\alpha_p \otimes \bar S^{\circ p}(b)\big ) = \bar \omega(b)  \EM_\omega^{(p)} (\alpha_p \otimes e )  + \EM_\omega^{(p)}\big (\alpha_p\otimes (\bar S^{\circ p}(b) - \bar \omega(b) e )\big ),
\end{equation}
and, by using the fact that $\EM_\omega^{(p)}$ are all contractions, we obtain
\begin{equation}\label{Eq:X103}
\begin{aligned}
& \|\EM_\omega^{(2p)} \big ( (\alpha_p \otimes 1^{\otimes p}) \otimes b \big ) - \bar \omega(b) b'\|_1^{\rm osy} \\
& \qquad \leq |\bar \omega(b)| \| \EM_\omega^{(p)} (\alpha_p \otimes e )-b'\|_1^{\rm osy}  +\| \alpha_p\| \|\bar S^{\circ p}(b) - \bar \omega(b) e )\|_1^{\rm osy}.
\end{aligned}
\end{equation}
Since the left side of the inequality goes to zero as $p \to \infty$, we can conclude that the left hand side of Eq.~\eqref{Eq:Full} contains the linear space generated by $b'$. Since $b'$ was arbitrary, the statement follows.\end{proof}

\subsection{Concluding remarks and a look ahead} The reduction and factorization processes described in this section produced the data consisting of the following:
\begin{enumerate}[\ (1)]
\item The local nuclear $C^\ast$-algebra $\Aa$.
\item The reduced space $\Bb_\omega$ with the structure of an operator system.
\item The u.c.p. map $\EM_\omega : \Aa \otimes \Bb_\omega \rightarrow \Bb_\omega$.
\item The u.c.p. functional $\bar \omega : \Bb_\omega \rightarrow \CM$.
\end{enumerate} 
It is certainly appropriate to say that the initial data $(\Aa_\ZM, \omega)$ was reduced and factorized to the data $(\Aa,\Bb_\omega,\EM_\omega,\bar \omega)$. 

We have also seen that questions related to the ergodicity of the state $\omega$ can be answered if the reduced data displays specific characteristics. Furthermore, our investigation of the ergodic states led us naturally to the concept of a full $\EM$ map and a connection was made between this concept and a condition that assures that the state is ergodic. This property will appear again, in an essential way, in the reconstruction phase of our program, investigated next. 

\section{Reconstruction Process}
\label{Sec:Reconstruction}

 In this section, we consider a set of data that shares the same attributes as the reduced data of a state over $\Aa_\ZM$. As we shall see, such data always produces a state on $\Aa_\ZM$ by a straightforward algorithm. However, the reduced data corresponding to the so-constructed state may not coincide with the initial data supplied as the input for the algorithm. The discussion at page 451 in \cite{FannesCMP1992} where $\Bb_\omega$ is taken as $\Aa_R$ supplies such an example. In such cases, we cannot draw any conclusions about the ergodicity of the state by simply examining the initial data. For this reason it is imperative to identify input data that leads to states that reduce back to the input data, which will then enable us to apply the results obtained in subsection~\ref{Sub:Asy}. These issues are explored to the fullest in this section. 

\subsection{Reconstruction algorithm}

Here we prove one of our main results.

\begin{theorem}\label{Th:Main1} Assume:
\begin{enumerate}[$\ \circ$]
\item A unital nuclear $C^\ast$-algebra $\Aa$;
\item An Archimedean matrix-ordered space $(\Ss,e)$;
\item A u.c.p. map $\EM : \Aa \otimes \Ss \rightarrow \Ss$;
\item A u.c.p. functional $\phi$ on $\Ss$.
\end{enumerate}
Then:
\begin{enumerate}[\rm i)]
\item Let $\EM^{(p)} : \Aa^{\otimes p} \otimes \Ss \rightarrow \Ss$ be the systems of maps defined iteratively as in Proposition~\ref{Pro:IsIter}, specifically,
\begin{equation}
\EM^{(1)} = \EM, \quad \EM^{(p+1)} = \EM\circ ({\rm id} \otimes \EM^{(p)}), \quad p \geq 1.
\end{equation}
Then the tower of linear functionals
\begin{equation}\label{Eq:OmegaTower}
\omega_{(p)}: \Aa^{\otimes p} \rightarrow \CM, \quad \omega_{(p)}=\phi \circ \EM^{(p)} \circ J_p, \quad p \geq 1,
\end{equation}
define a state $\omega_R$ on $\Aa^{\otimes \NM^\times} \simeq \Aa_R$, where $J_p$'s are the unital complete order embeddings $J_p:\Aa^{\otimes p} \rightarrowtail \Aa^{\otimes p} \otimes \Bb_\omega$.

\vspace{0.1cm}

\item Suppose $\phi \circ \bar S = \phi$, where $\bar S$ is defined as in Proposition~\ref{Pro:Cluster2}, specifically,
\begin{equation}
\bar S:\Ss \rightarrow \Ss, \quad \bar S = \EM \circ \bar L,
\end{equation} 
with  $\bar L$ being the unital complete order embedding $\bar L : \Ss \rightarrow \Aa \otimes \Ss$. Then there exists a unique shift-invariant state $\omega$ on $\Aa_\ZM$ extending $\omega_R$ derived at point {\rm i)}.
\end{enumerate}
\end{theorem} 

\begin{proof} i) As compositions of u.c.p. maps (recall point T3 of Definition~\ref{Def:TP}), the maps $\EM^{(p)}$ are u.c.p. and the maps $J_p$ are as well u.c.p. As such, $\omega_{(p)}$ are u.c.p., hence states on $\Aa^{\otimes p}$. Furthermore, if $\mathfrak j_{q,p} : \Aa^{\otimes p} \rightarrowtail \Aa^{\otimes p} \otimes \Aa^{\otimes (q-p)} = \Aa^{\otimes q}$ are the standard $C^\ast$-algebra embeddings for $q \geq p$, then $\omega_{(q)} \circ \mathfrak j_{q, p} = \omega_{(p)}$ for any $q \geq p$ and, as such, the tower of states respects the structure maps of the directed tower of $C^\ast$-algebras and, as such, it defines a state on the inductive limit, which is $\Aa^{\otimes \NM^\times} \simeq \Aa_R$. 

ii) We denote by $\omega_R$ the inductive limit of states from point i), and let $S_R : \Aa^{\otimes \NM^\times} \to \Aa^{\otimes \NM^\times}$ be the $C^\ast$-algebra morphism $a_R \mapsto 1 \otimes a_R$. Then, for any $p \geq 1$ and $\alpha \in \Aa^{\otimes p}$, we have
\begin{equation}
(\omega_R \circ S_R)\big (\mathfrak j_{(\infty,p)}(\alpha)\big ) = \omega_{(p+1)} \big (1 \otimes \alpha \big) = \big( \phi \circ \EM^{(p+1)}\big)\big ( ( 1 \otimes \alpha ) \otimes e \big ).
\end{equation}
Using the very definition of $\EM^{(p)}$, we find
\begin{equation}
(\omega_R \circ S_R)\big (\mathfrak j_{(\infty,p)}(\alpha)\big ) = (\phi \circ \EM)\big (1 \otimes \EM^{(p)}(\alpha \otimes e )\big),
\end{equation}
and, after invoking the definition of $\bar L$,
\begin{equation}
(\omega_R \circ S_R)\big (\mathfrak j_{(\infty,p)}(\alpha)\big ) = (\phi \circ \EM\circ \bar L)\big (\EM^{(p)} (\alpha \otimes e)\big).
\end{equation}
Since $\EM \circ \bar L = \bar S$ and $\phi \circ \bar S = \phi$, we must have $\omega_R \circ S_R = \omega_R$. Any state with such property can be uniquely extended to a shift-invariant state over $\Aa_\ZM$.
\end{proof}

Let $(\Aa,\Bb_\omega,\EM_\omega,\bar \omega)$ be the reduced data for the state produced via Theorem~\ref{Th:Main1} from the input $(\Aa,\Ss,\EM,\phi)$ with $\phi$ shift-invariant. As already stated, without additional information about the input data, there is no reason to assume that the two sets of data coincide. It is imperative to find out when the two actually coincide. A related problem is the identification of states that can be produced with the algorithm from Theorem~\ref{Th:Main1}, for some input data. The two mentioned issues are actually related: 

\begin{proposition}\label{Pro:Unique} If $\omega$ and $\omega'$ are shift invariant states and the data $(\Aa_\ZM, \omega)$ and $(\Aa_\ZM,\omega')$ both reduce to the same data $(\Aa,\Bb,\EM,\phi)$, then necessarily $\omega = \omega'$.
\end{proposition}

\begin{proof} Let $\omega_R$ be the restriction of $\omega$ on $\Aa_R$. From Proposition~\ref{Pro:Omega}, we know that $\omega_R = \phi \circ q$ and, from Proposition~\ref{Pro:IndL}, we know that $q$ is completely determined by $\EM$. Since the same arguments apply for $\omega'$, the conclusion is that $\omega_R$ and $\omega'_R$ coincide on a dense subspace of $\Aa_R$, hence on the whole $\Aa_R$. Due their shift-invariance, they must also coincide on the whole $\Aa_\ZM$.\end{proof}

\begin{corollary} Let $\omega$ be a state over $\Aa_\ZM$ and $(\Aa,\Bb_\omega,\EM_\omega,\bar \omega)$ be its reduced data. Then the algorithm from Theorem~\ref{Th:Main1} with input $(\Aa,\Bb_\omega,\EM_\omega,\bar \omega)$ generates back the state $\omega$.
\end{corollary}

The following statement supplies sufficient conditions for an input data set $(\Aa,\Ss,\EM,\phi)$ to be the reduced data of some state over $\Aa_\ZM$

\begin{theorem}\label{Th:Main2} Assume the conditions of Theorem~\ref{Th:Main1} and, additionally, that $\phi$ is shift invariant and $\EM$ is full,
\begin{equation}\label{Eq:FullE}
\overline{\bigcup\nolimits_{p \geq 1} \bigcup\nolimits_{x \in \Aa^{\otimes p}} \EM^{(p)}(x \otimes s)} = \Ss, \quad \forall \ s \in \Ss, \ s \neq 0.
\end{equation} 
Let $\omega$ be the shift-invariant state produced by the algorithm from Theorem~\ref{Th:Main1} from the input data $(\Aa,\Ss,\EM,\phi)$. Then the data $(\Aa_\ZM,\omega)$ reduces back to the input data $(\Aa,\Ss,\EM,\phi)$. In particular, this is the case if $\bar S$ satisfies relation~\eqref{Eq:FixP}.
\end{theorem}

\begin{proof}
We first note that $(\EM^{(q)} \circ J_q) \circ \mathfrak j_{q, p} = \EM^{(p)} \circ J_p$ for any $q \geq p$ and, as such, we have a direct limit map $\EM^{(\infty)} \circ J_\infty : \Aa^{\otimes \NM^\times} \to \Ss$. Note that $\EM^{(\infty)} \circ J_\infty$ is unital and contractive and continues to enjoy Markov's property~\eqref{Eq:Markov}. We define the closed subset 
\begin{equation}
\Kk = {\rm Ker}\big (\EM^{(\infty)} \circ J_\infty\big ) \subset \Aa^{\otimes \NM^\times}.
\end{equation}
We will show that $\Kk$ coincides with the entanglement kernel $\Kk_\omega$ of the state $\omega$. For $x = \beta_{(-p+1,0)}(\alpha)$ with $\alpha \in \Aa^{\otimes p}$, we have
\begin{equation}\label{Eq:X23}
\omega(x \cdot a_R)=\phi \Big (\EM^{(p)} \big (\alpha \otimes (\EM_\infty \circ J_\infty) (a_R) \big ) \Big ), \quad \forall \ a_R \in \Aa^{\otimes \NM^\times} \simeq \Aa_R.
\end{equation}
As such, if $a_R \in \Kk$, then $\omega_x(a_R)=0$ for $x$ in a dense subset of $\Aa_L$, hence for all $x$ in $\Aa_L$ because of the continuity of the state. This proves that $\Kk \subseteq \Kk_\omega$. The reciprocal is also true in the stated conditions. Indeed, if $a_R \in \Kk_\omega$, then necessarily
\begin{equation}\label{Eq:X73}
\phi \Big (\EM^{(p)} \big (\alpha \otimes (\EM_\infty \circ J_\infty) (a_R) \big ) \Big )=0
\end{equation}
for any $\alpha \in \Aa^{\otimes p}$ and $p \geq 1$, which is a direct consequence of the identity~\eqref{Eq:X23}. Let $s=(\EM_\infty \circ J_\infty) (a_R) \in \Ss $ and assume that $s \neq 0$. In this case, since $\EM$ is full,
\begin{equation}\label{Eq:Surj1}
\overline{\bigcup\nolimits_{p \geq 1} \bigcup\nolimits_{\alpha \in \Aa^{\otimes p}} \EM^{(p)}(\alpha \otimes s)} = \Ss,
\end{equation}
and \eqref{Eq:X73} can be true only if $\phi =0$. This contradiction proves that $s=0$ or, in other words, that $a_R \in \Kk$.

We have established that the entanglement kernel $\Kk_\omega$ of $\omega$ coincide with the kernel ${\rm Ker}(\EM^{(\infty)} \circ J_\infty)$. Then the projection $q : \Aa_R \to \Aa_R/K_\omega$ coincides with $\EM^{(\infty)} \circ J_\infty$. Our last task is to show that $\EM$ given at the start of the reconstruction process coincides with the map defined in Eq.~\eqref{Eq:EDef} for the reconstructed state $\omega$. Specifically, we need prove that
\begin{equation}\label{Eq:TLast}
\EM\big(a \otimes (\EM^{(\infty)} \circ J_\infty)(a_R) \big) = (\EM^{(\infty)} \circ J_\infty) \big (S(\beta_{(0,0)}(a)a_R)\big ),
\end{equation} 
for all $a \in \Aa$ and $a_R \in \Aa^{\otimes \NM^\times}$. For any $\alpha \in \Aa^{\otimes p}$, we have
\begin{equation}
\begin{aligned}
\EM \Big (a \otimes \big(\EM^{(\infty)} \circ J_\infty\big) \big( \mathfrak j_{\infty,p}(\alpha)\big) \Big )
& =\EM \Big (a \otimes \big(\EM^{(p)} \circ J_p\big)(\alpha) \Big ) \\
& = \EM^{(p+1)}\Big ( \big (a \otimes \alpha \big ) \otimes e \big).
\end{aligned}
\end{equation}
The last line can be cast as $(\EM^{(\infty)} \circ J_{\infty}) \big (\mathfrak j_{\infty,p+1} (a \otimes \alpha) \big)$ and $j_{\infty,p+1} (a \otimes \alpha)$ can be seen to coincide with $S(\beta_{(0,0)}(a)a_R)$, under the isomorphism $\Aa^{\otimes \NM^\times} \simeq \Aa_R$. As such, the relation~\eqref{Eq:TLast} holds for a dense subset of $\Ss$, hence on the whole $\Ss$.
\end{proof}

\begin{example}{\rm It is easy to exemplify how the statement in Theorem~\ref{Th:Main2} fails if the full-ness condition is not satisfied. Indeed, let $\EM : \Aa \otimes \Ss \to \Ss$ be full as stated. Consider 
\begin{equation}
\hat \EM : \Aa \otimes \Ss \otimes \Bb \to \Ss \otimes \Bb, \quad \hat \EM = \EM \otimes {\rm id},
\end{equation} 
where $\Bb$ is an auxiliar operator system. A short calculation shows that 
\begin{equation}
\hat \EM^{(p)}(\alpha \otimes s \otimes b) = \EM^{(p)}(\alpha \otimes s) \otimes b,
\end{equation}
hence it is evident that $\hat \EM$ fails to be full. Now, consider the state $\hat \omega=\bar \omega \otimes \xi$ over $\Ss \otimes \Bb$, with $\xi$ any state on $\Bb$. Then $(\hat \EM,\hat \omega)$ produces the same state $\omega$ over $\Aa_\ZM$ through the reconstruction process, but $(\Aa_\ZM,\omega)$ reduces to $\Ss$ and not to $\Ss \otimes \Bb$.
}$\Diamond$
\end{example}

\begin{corollary} Assume the conditions of Theorem~\ref{Th:Main1} and, additionally, that $\phi \circ \bar S = \phi$, that the map $\EM$ is full and that the input data displays the asymptotic clustering property
 \begin{equation}
\lim_{r \rightarrow \infty}\sup \Big \{ \Big | \phi \Big(\EM^{(p)}\big (\alpha \otimes \bar S^{\circ r}(s - \phi(s) \, e) \big ) \Big )\Big | \Big \}=0, \quad \forall \ s \in \Ss,
\end{equation}
where the supremum is over all $n \geq 1$, $\alpha \in \Aa^{\otimes p}$ with $\|\alpha\|=1$. Then the reconstructed state $\omega$ displays the asymptotic clustering property and, as such, it is ergodic.
\end{corollary}

\begin{proof} In the stated conditions, the reduced data $(\Aa,\Bb_\omega,\EM_\omega,\bar \omega)$ for reconstructed state $\omega$ coincides with the input data $(\Aa,\Ss,\EM,\phi)$. Then we can apply the results of subsection~\ref{Sub:Asy} and the statement follows from Proposition~\ref{Pro:Cluster2}.\end{proof}

We now put forward the ideal scenario that will sought for in all our examples:

\begin{theorem}\label{Th:Main3} Let $(\Aa,\Ss,\EM,\phi)$ be as in Theorem~\ref{Th:Main1} and suppose that
\begin{equation}
\lim_{r \to \infty} \bar S^{\circ r}(s) \to \phi(s) \, e, \ \forall \ s \in \Ss,
\end{equation}
checks for the input data. Then:
\begin{enumerate}[{\rm i)}]
\item $\phi$ is shift invariant, $\phi \circ \bar S = \phi$;
\item The map $\EM$ is full and, as such, the data $(\Aa_\ZM,\omega)$, with $\omega$ the shift-invariant state generated by the algorithm form Theorem~\ref{Th:Main1}, reduces back to the input data $(\Aa,\Ss,\EM,\phi)$;
\item The state $\omega$ displays the asymptotic clustering property and, as such, it is ergodic.
\end{enumerate}
\end{theorem}

\begin{proof} i) is evident. ii) follows from the same argument as in Proposition~\ref{Pro:EDense}. iii) follows from Proposition~\ref{Pro:Cluster2}, which applies to the present context due to ii).\end{proof} 

\subsection{Examples of reconstructed states}\label{Sec:Reco}

Below, we give examples of ergodic states derived with the algorithm described in Theorem~\ref{Th:Main1}.

\begin{example}\label{Ex:ProdState3}{\rm Let us consider the product state discussed in Examples~\ref{Ex:ProdState1} and \ref{Ex:ProdState2}, where we found $\Bb_\omega = \CM \cdot e$ and $\EM_\omega(a \otimes e) = \omega_0(a) \cdot e$. Hence, $\EM_\omega$ is obviously full and there is only one state on $\Bb_\omega$, $\bar \omega(e) =1$. As such, 
\begin{equation}
\omega_{(1)}(a) : = \bar \omega \big(\EM(a \otimes e)\big) = \omega_0(a)
\end{equation}
and
\begin{equation}
\omega_{(2)}(a_1 \otimes a_2) : = \bar \omega\big(\EM^{(2)}(a_1 \otimes a_2 \otimes e)\big) = \bar \omega\big(\EM\big(a_1 \otimes \EM(a_2)\big )\big) = \omega_0(a_1) \omega_0(a_2).
\end{equation}
Iterating further, one finds
\begin{equation}
\omega_{(n)}(a_1 \otimes \cdots \otimes a_n) = \omega_0(a_1) \cdots \omega_0(a_n),
\end{equation}
which confirms that the product state $\omega_0^{\otimes \ZM}$ is indeed reproduced by the reconstruction algorithm.
}$\Diamond$
\end{example}

\begin{example}\label{Ex:AKLT2}{\rm We reconsider here the class of AKLT states introduced in Example~\ref{Ex:AKLT1}, which we now can analyze fully, without assuming any convergence of periodic approximants. We recall that the setting was that of a nuclear $C^\ast$-algebra $\tilde \Aa$, of a projection $p$ from $\tilde \Aa \otimes \tilde \Aa$, and of a a positive map $\xi_0 : \tilde \Aa \otimes \tilde \Aa \rightarrow \CM$. Both operator systems $\Aa$ and $\Ss$ are defined in terms of this data, namely $\Aa = p( \tilde \Aa \otimes \tilde \Aa)p$, with $p$ standing for the unit of $\Aa$, and $\Ss = \tilde \Aa$. We let $\mathfrak j : p \hat \Aa p \rightarrow \hat \Aa$ be the {\it non-unital} embedding that takes $p$ into $p$ rather than into the unit of $\hat \Aa$. As a $C^\ast$-algebra morphism, $\mathfrak j$ is a c.p. map. Lastly, we define $\EM$ as
\begin{equation}\label{Eq:EDeF}
\EM : p ( \tilde \Aa \otimes \tilde \Aa) p \otimes \tilde \Aa \rightarrow \Aa, \quad \EM =  ({\rm id} \otimes \xi_0) \circ (\mathfrak j \otimes {\rm id}),
\end{equation}
which, as a composition of two c.p. maps,\footnote{Property T3 in Definition~\ref{Def:TP} assures us that this is indeed the case.} is a c.p. map. This map is also unital, provided 
\begin{equation}
({\rm id} \otimes \xi_0)(p \, \otimes \tilde 1) = \tilde 1.
\end{equation} 
Therefore, whenever $p$ and $\xi_0$ fulfill this constraint, an AKLT-type state can be reconstructed from the data 
\begin{equation}
{\rm AKLT}= \Big (\Aa := p(\tilde \Aa \otimes \tilde \Aa)p, \ \Ss : = \tilde \Aa, \ \EM : =({\rm id} \otimes \xi_0) \circ (\mathfrak j \otimes {\rm id}), \ \phi \Big ),
\end{equation} 
where $\phi$ is a shift invariant state on $\tilde \Aa$. This machinery now works equally for finite and infinite dimensional nuclear $C^\ast$-algebras $\Aa$. Below, we give two examples where the AKLT data satisfies the conditions of Theorem~\ref{Th:Main3}.
}$\Diamond$
\end{example}

\begin{remark}{\rm The particular case studied in \cite{AffleckCMP1988}, illustrated in Fig.~\ref{Fig:1DSystem} and partially analyzed in Example~\ref{Ex:AKLT1}, corresponds to $\tilde \Aa = M_2(\CM)$, hence $\hat \Aa = M_2(\CM) \otimes M_2(\CM) \simeq M_4(\CM)$, and $p$ is the rank-3 projection
\begin{equation}\label{Eq:PP1}
p = \tfrac{3}{4} \, \sigma_0 \otimes \sigma_0 + \tfrac{1}{4} \sum_{i=1}^3 \sigma^i \otimes \sigma ^i \in M_4(\CM),
\end{equation}
such that $\Aa = p M_4(\CM) p \simeq M_3(\CM)$. Above, $\sigma_i$, $i=1,2,3$, are Pauli's matrices and $\sigma_0$ is the identity of $M_2(\CM)$. Furthermore,
\begin{equation}
\xi_0(m) = \tfrac{4}{3} \, {\rm Tr}\big ( m p_0 \big ), \quad m \in M_4(\CM),
\end{equation}
with $p_0$ being the rank-1 projection
\begin{equation}\label{Eq:P0}
p_0=\tfrac{1}{4} \, \sigma_0 \otimes \sigma_0 - \tfrac{1}{4} \sum_{i=1}^3 \sigma^i \otimes \sigma ^i \in M_4(\CM).
\end{equation}
In this case, \eqref{Eq:EDeF} takes the following concrete form:
\begin{equation}\label{Eq:EMSigma}
\EM\Big (p\big (\gamma \otimes \gamma' \big ) p \otimes \gamma'' \Big ) = \tfrac{4}{3}\sum_{\alpha,\beta=0}^3 g_{\alpha} g_{\beta}  {\rm Tr} \Big ( \big ( \sigma^\alpha \gamma' \sigma^{\beta} \otimes \gamma'' \big ) p_0 \Big ) \, \sigma^{\alpha} \gamma \sigma^{\beta},
\end{equation}
where $g_{0}=\tfrac{3}{4}$ and $g_{i}=\tfrac{1}{4}$ for $i=\overline{1,3}$ and $\gamma,\gamma',\gamma'' \in M_2(\CM)$. We then find
\begin{equation}
\EM(p \otimes \sigma_0) = \sigma_0, \quad \EM(p \otimes \vec \alpha \cdot \vec \sigma) = -\tfrac{1}{3} \, \vec \alpha \cdot \vec \sigma,
\end{equation}
where $\alpha \in \CM^3$ and $\vec \alpha \cdot \vec \sigma = \sum_{i=1}^3 \alpha_i \, \sigma_i$. Hence, the data satisfy the constraint mentioned in Example~\ref{Ex:AKLT2}. Furthermore, since any $s \in M_2(\CM)$ can be written uniquely as $s = \alpha_0 \, \sigma_0 + \vec \alpha \cdot \vec \sigma$, we have
\begin{equation}
\bar S^{\circ r}( s) = \EM(p \otimes \EM( \ldots \EM(p \otimes s)\ldots )) = \alpha_0 \, \sigma_0 + \big ( -\tfrac{1}{3} \big )^r \, \vec \alpha \cdot \vec \sigma,
\end{equation}
hence $\bar S$ satisfies relation~\eqref{Eq:FixP}. As a consequence, there is only one $\bar S$-invariant functional $\bar \omega$ on $M_2(\CM)$, which is
\begin{equation}
\bar \omega(\alpha_0 \, \sigma_0 + \vec \alpha \cdot \vec \sigma) = \alpha_0,
\end{equation}
and the reconstructed state $\omega$ is ergodic. Two-site correlation functions can be also computed explicitly, by following Example~\ref{Re:CorrF}:
\begin{equation}
C(a,a';r) = (\bar \omega\circ \EM)\big (a \otimes \bar S^{\circ r}(\hat a') \big)=\omega(a)\omega(a') + \big ( -\tfrac{1}{3} \big )^r (\bar \omega \circ \EM)(a, \vec \alpha \cdot \vec \sigma),
\end{equation}
where $\vec \alpha \cdot \vec \sigma = \EM(a' \otimes \sigma_0) - \omega(a')$ and $\omega(a) = (\bar \omega \circ \EM)(a \otimes \sigma_0)$. For given specific entries $a$ and $a'$, the calculation can be completed by applying~\eqref{Eq:EMSigma}.
}$\Diamond$
\end{remark}

\begin{example}{\rm We can take a page from the above example and produce an AKLT data with $\Aa$ infinite dimensional. For this, we want to engage the full $C^\ast$-algebra $C^\ast G$ of an infinite amenable discrete group $G$,\footnote{For amenable groups, the full and reduced group $C^\ast$-algebras coincide and are nuclear \cite{LanceJFA1973}.} but, typically, such algebras are poor in projections, unless we tensor them with a matrix algebra. For this reason, the AKLT construction requires a slight generalization, which we now explain. We start with a projection $\tilde p \in  M_N(\CM) \otimes C^\ast G$. The group morphism $G \ni g \mapsto g \times g \in G \times G$ lifts canonically to an algebra morphism $M_N(\CM) \otimes C^\ast G \to M_N(\CM) \otimes ( C^\ast G \otimes C^\ast G)$, hence, the image of $\tilde p$ through this map is still a projection, which we take as the projection $p$ in the AKLT construction. The latter can be presented as a norm convergent series 
\begin{equation}
p=\sum_{g \in G} c_g \cdot g \otimes g, \ c_g \in M_N(\CM), \ c_{g^{-1}} = c_g^\ast.
\end{equation} 
We now can specify half of the input data:
\begin{equation}
\Aa = p(M_N(\CM) \otimes C^\ast G \otimes C^\ast G)p, \quad \Ss = C^\ast G.
\end{equation}
Next, we choose a positive element $\tilde q_0$ from $C^\ast G$ with $\|q_0\| \leq 1$, and push it through the same morphism into  $M_N(\CM) \otimes ( C^\ast G \otimes C^\ast G)$. If 
\begin{equation}
q_0=\sum\nolimits_{g \in G} d_g \cdot g \otimes g, \ d_g \in \CM, \ d_{g^{-1}}=d_g^\ast,
\end{equation} 
is the result of that action, then we define the positive map
\begin{equation}
C^\ast G \otimes C^\ast G \ni \sigma \mapsto \xi_0(\sigma) = (\Tt \otimes \Tt)(\sqrt{q_0} \, \sigma  \sqrt{q_0})=(\Tt \otimes \Tt)(\sigma   q_0) \in \CM,
\end{equation}
where $\Tt$ is the standard trace on $C^\ast G$, $\Tt(\sum_g \lambda_g \cdot g) = \lambda_e$. Lastly, we define
\begin{equation}
\EM : \Aa \otimes \Ss \to \Ss, \quad \EM = ({\rm tr} \otimes {\rm id} \otimes \xi_0)\circ (\mathfrak j \circ {\rm id}),
\end{equation}
where ${\rm tr}$ is the trace state on $M_N(\CM)$ and $\mathfrak j$ is the non-unital $C^\ast$-algebra embedding $\Aa \rightarrowtail M_N(\CM) \otimes C^\ast G \otimes C^\ast G$. As a composition of c.p. maps, $\EM$ is a c.p. map. Furthermore, if $s=\sum_{g \in G} s_g \cdot g \in \Ss$, then we have
\begin{equation}
\begin{aligned}
({\rm tr} \otimes {\rm id} \otimes \xi_0)(p \otimes s) &= ({\rm tr} \otimes {\rm id} \otimes \Tt \otimes \Tt)\big (\sum_{f,g,h\in G} c_g s_h d_f \cdot g \otimes gf \otimes hf\big )
\end{aligned}
\end{equation}
and, by using the fact that $\Tt(g) = \delta_{g,e}$, we find
\begin{equation}
\bar S(s) = \EM(p \otimes s) = \sum\nolimits_{g \in G}   d_{g}^\ast \, {\rm tr}(c_g)s_g \cdot g.
\end{equation}
Therefore, if we choose the coefficients of $q_0$ such that
\begin{equation}
d_e = 1/{\rm tr}(c_e), \quad |d_{g}| < 1/|{\rm tr}(c_g)| \ {\rm for} \ g \neq e,
\end{equation}
which we can always do, then $\EM(p \otimes e) = e$ and
\begin{equation}
\bar S^{\circ r}(s)  = \sum\nolimits_{g \in G}   \big (d_{g}^\ast {\rm tr}(c_g)\big )^rs_g \cdot g \to s_e \cdot e, \ {\rm as} \ r \to \infty, \quad \forall \ s \in \Ss.
\end{equation}
The conclusion is that the reconstruction algorithm for the modified AKLT data
\begin{equation}
{\rm mAKLT}= \Big (\Aa := p(M_N(\CM) \otimes \tilde \Aa \otimes \tilde \Aa)p, \ \Ss : = \tilde \Aa, \ \EM : =({\rm tr} \otimes {\rm id} \otimes \xi_0) \circ (\mathfrak j \otimes {\rm id}) \Big ),
\end{equation} 
produces an ergodic state $\omega$ on $\Aa^{\otimes \ZM}$.
}$\Diamond$
\end{example}

\begin{remark}{\rm The last example can be given the physical interpretation of the thermodynamic limit of states for $N$ distinguishable quantum particles hopping on the Cayley graph of $G$, where the states are invariant against circular permutation of the identities of the particles.
}$\Diamond$
\end{remark}

\section{Stinespring Representations of Operator Product States}\label{Sec:Stinespring}

The goal of this section is to appeal to the Stinespring representation of map $\EM$ and derive a possibly simpler set of input data for a reconstruction algorithm (see Proposition~\ref{Pro:TEM}). It is argued that this new type of input data samples densely the space of input data introduced in the previous section, provided $\Ss$ can be embedded in a postliminal $C^\ast$-algebra.

\subsection{Background: Types of $C^\ast$-algebras and their spectra} When appealing to the Stinespring representation of the map $\EM$, one immediately encounters the space of representations of the $C^\ast$-algebras involved. As we shall see, the generic representation of $\EM$ engages many if not all irreducible representations of an embedding $C^\ast$-algebra for $\Bb_\omega$. This calls for a brief overview of the spaces where the representations of $C^\ast$-algebras live.

A representation of a $C^\ast$-algebra $\Qq$ is a $C^\ast$-algebra morphism $\pi: \Qq \to B(H)$ for some Hilbert space $H$. Both the $C^\ast$-algebras and the Hilbert spaces will be assumed separable in this section. When referring to a specific representation $\pi$, we will use the symbol $H_\pi$ for the associated Hilbert space. We recall that two representations are called equivalent if there exists a unitary map between the corresponding Hilbert spaces intertwining the two representations.

\begin{definition} A representation $\pi$ of a $C^\ast$-algebra $\Qq$ is called topologically irreducible if $0$ and $H_\pi$ are the only closed subspaces invariant for the action of the algebra $\Qq$. The set of equivalence classes of irreducible representations defines the spectrum of the algebra, typically indicated by a hat as in $\hat \Qq$.
\end{definition}

The spectrum of a $C^\ast$-algebra accepts a canonical topology, which can be introduced in a multitude of distinct ways (see \cite{DixmierBook1}~Ch.~3, \cite{FellDoranBook}~Ch. VII, \cite{RaeburnWilliamsBook1998}~Ch.~ A). We want to describe this topology here.

\begin{definition}[\cite{DixmierBook1},~Sec.~3.1] A closed double-sided ideal of a $C^\ast$-algebra $\Qq$ is called primitive if it is the kernel of an irreducible representation of $\Qq$. The set of primitive ideals is usually specified as ${\rm Prim}(\Qq)$. 
\end{definition}

\begin{proposition}[\cite{DixmierBook1},~Sec.~3.1]\label{Prop:JacobTop}
For each subset $T \subseteq {\rm Prim}(\Qq)$, one lets $I(T)$ be the intersection of the elements of $T$ and considers the set $\overline T$ of all primitive ideals of $\Qq$ containing $I(T)$. Then there exists a unique topology on ${\rm Prim}(\Qq)$ such that $\overline T$ is the topological closure of $T$. 
\end{proposition}

\begin{definition} The topology defined by Proposition~\ref{Prop:JacobTop} is called the Jacobson topology. The ${\rm Prim}$ space is endowed with this natural topology.
\end{definition}

\begin{remark}{\rm The naturality stems from the observation that the closed subsets of ${\rm Prim}(\Qq)$ and the closed double-sided ideals of $\Aa$ are in a bijective relation, established by the map ${\rm Prim}(\Qq) \ni T = \overline T \mapsto \bigcap \{J, \ J \in T\}$. In other words, the lattice of ideals of $\Qq$ can be derived from the Jacobson topology of $Prim(\Qq)$.
}$\Diamond$
\end{remark}

\begin{definition}\label{Def:HatTopo}  Since the kernels of equivalent representations coincide, there exist a natural surjective map $\hat \Qq \ni [\pi] \mapsto {\rm ker}\, \pi \in {\rm Prim}(\Qq)$ and the topology of $\hat \Qq$ is defined to be the pull-back topology through this map.
\end{definition}

The pure states over a $C^\ast$-algebra $\Qq$ supply irreducible representations. If $P(\Qq)$ stands for the set of pure states, then this set comes equipped with the weak-$\ast$ topology inherited from the dual of $\Qq$.

\begin{proposition}[\cite{DixmierBook1}~3.4.12] The topology of $\hat \Qq$ introduced in Definition~\ref{Def:HatTopo} coincides with the quotient topology of the topology of $P(\Qq)$ for the canonical surjective map $P(\Qq) \to \hat \Qq$.
\end{proposition} 

A topological space is called a $T_0$-space if for any two distinct points of the space there is a neighborhood of one of the points which does not contain the other. While Prim spaces are always $T_0$ \cite{DixmierBook1}[Prop.~3.1.3], this is not always the case for the spectra. When this does happen, we have: 

\begin{proposition}[\cite{DixmierBook1},~Prop.~3.1.6]\label{Pro:T0} The following three conditions are equivalent:
\begin{enumerate}[\rm (i)]
\item $\hat \Qq$ is a $T_0$-space.
\item Two irreducible representations of $\Qq$ with the same kernel are equivalent.
\item The canonical map $\hat \Qq \to {\rm Prim}(\Qq)$ is a homeomorphism.
\end{enumerate}
\end{proposition}

The simplest spectrum is that consisting of a single point. The separable $C^\ast$-algebras displaying such spectra are precisely the elementary ones, that is, the ones that are isomorphic to the algebra of compact operators over some Hilbert space (finite or infinite). These were exactly the $C^\ast$-algebras engaged in the study of finitely-correlated states in \cite{FannesCMP1992}, and they remain very relevant for the more general context considered here (see Example~\ref{Ex:HilbertC}). All finite $C^\ast$-algebras have finite Hausdorff spectra \cite{BunceTAMS1975}[Corollary~8]. Next up in the ladder of complexity come the dual $C^\ast$-algebras \cite{KaplanskyAM1948} (see also \cite{DixmierBook1}[4.7.20] for a brisk characterization), which all have discrete spectra. The $C^\ast$-algebras of compact groups are dual \cite{KaplanskyAM1948} and, in fact, this property characterizes completely the compact groups among the locally compact groups \cite{ElliotPAMS1971,BaggettJFA1972}. Liminal and postliminal $C^\ast$-algebras \cite{DixmierBook1}[Ch.~4] also have a fairly complete characterization of their spectra. In particular, postliminal $C^\ast$-algebras have only type I representations and satisfy the equivalent properties stated in Proposition~\ref{Pro:T0} \cite{GlimmAM1961}. Furthermore:

\begin{proposition}[\cite{DixmierBook1},~Ch.~4]\label{Pro:Postliminal} Let $\Qq$ be a separable postliminal $C^\ast$-algebra. Let $\pi$ be a non-degenerate representation of $\Qq$ on a separable Hilbert space. Then there exist mutually singular positive measures $\bm \mu_1$, $\bm \mu_2$, \ldots, $\bm \mu_\infty$ on $\hat \Qq$, such that
\begin{equation}\label{Eq:Disinter}
\pi \simeq \int\nolimits^\oplus_{\hat \Qq} d\bm \mu_1(\xi) \, \xi \oplus 2\int\nolimits^\oplus_{\hat \Qq} d\bm \mu_2(\xi) \, \xi \oplus \cdots \oplus \aleph_0 \int\nolimits^\oplus_{\hat \Qq} d\bm \mu_\infty(\xi) \, \xi,
\end{equation}
where the coefficients in front of the integral signs indicate the multiplicities. The system of measures $\{\bm \mu_i\}$ is fixed by $\pi$ up to measure equivalence.
\end{proposition}

\begin{remark}{\rm In general, the integration and disintegration of representations of a $C^\ast$-algebra is developed over the quasi-spectrum of the algebra  equipped with its natural Mackey-Borel structure. For a separable postliminal $C^\ast$-algebra, however, the spectrum endowed with the topological Borel structure coincides with the Mackey-Borel structure
\cite{DixmierTAMS1962} (see also \cite{DixmierBook1}[Ch.~7]).
}$\Diamond$
\end{remark}

All commutative $C^\ast$-algebras are postliminal. $C^\ast$-algebras of type-I topological groups are also sources of postliminal $C^\ast$-algebras (see \cite{FollandBook2000}[Theorem~7.8] and \cite{DeitmarEchterhoffBook}[Example~8.5.1] for explicit lists of type-I groups). In particular, a discrete group is of type-I if and only if it contains an abelian subgroup of finite index. In particular, the $C^\ast$-algebra of the space-groups engaged in crystallography (i.e. the lattices of the Euclidean group) are postliminal.

We will restrict our discussion of the operator-product presentation of a state over $\Aa_\ZM$ to the situations where the operator system $\Bb_\omega$ can be embedded in a postliminal $C^\ast$-algebra. Besides all the above nice features, the class of postliminal $C^\ast$-algebras has the following special property:  

\begin{theorem}[\cite{RaeburnWilliamsBook1998}, Th.~B.45]\label{Th:TIrrep} Let $\Aa$ and $\Bb$ be nuclear $C^\ast$-algebras. Then the map $(\pi,\eta) \mapsto \pi \otimes \eta$ induces a homeomorphism of $\hat \Aa \times \hat \Bb$ onto its range in $\widehat {\Aa \otimes \Bb}$. Furthermore, if either $\Aa$ or $\Bb$ is postliminal, then this homeomorphism is surjective.
\end{theorem}

There are alternative tools to the disintegration theory summarized in Proposition~\ref{Pro:Postliminal}. For example, Fell and Doran describe in section~VI.14 of \cite{FellDoranBook} the notion of discretely decomposable $\ast$-representations. Such a representation accepts an essentially unique direct decomposition in irreducible representations. Another useful and related concept is that of approximately equivalent representations.

\begin{definition}[\cite{DavidsonBook},~pg.~57] Two representations $\pi$ and $\eta$ of a $C^\ast$-algebra $\Qq$ are said to be approximately unitary equivalent, written as $\pi \sim_a \eta$, if there exists a sequence of unitary transformations $U_n : H_\pi \to H_\eta$ such that
\begin{equation}\label{Eq:AppEquiv}
\pi(q) = \lim_{n \to \infty} U_n^\ast \, \eta(q) \, U_n, \quad \forall \ q\in \Qq,
\end{equation}
where the limit is in the norm topology of $B(H_\pi)$. Note that the difference $\pi(q) - U_n^\ast \eta(q) U_n$ can be arranged to be in $K(H_\pi)$ if $\pi \sim_a \eta$.
\end{definition}

Approximately unitary equivalence fits our purposes because of the following:

\begin{proposition}[\cite{DavidsonBook}~Corollary~II.5.9]\label{Pro:AUEq} Every representation of a separable $C^\ast$-algebra on a separable Hilbert space is approximately unitarily equivalent to a direct countable sum of irreducible representations.
\end{proposition}

\begin{proposition} Let $\Phi: \Qq \to B(H)$ be a u.c.p. map and $\Phi(q) = V^\ast \pi(q) V$ be a Stinespring representation of $\Phi$ as in Theorem~\ref{Th:Stinespring}.\footnote{$H_\pi$ can be chosen separable when both $\Qq$ and $H$ are separable \cite{PaulsenBook2002}[pg.~45].} Let $U_n : H_\pi \to \sum^\oplus_{[\zeta]} H_\zeta$ be the unitary transformations implementing the approximate unitary equivalence mentioned in Proposition~\ref{Pro:AUEq}, where the direct sum runs over a countable subset of $\hat \Qq$. Then
\begin{equation}
\Phi_n : \Aa \to B(H), \quad \Phi_n(q) : = V^\ast_n \Big (\sum\nolimits^\oplus_{[\zeta]} \zeta(q) \Big ) V_n, \quad V_n : = U_n V,
\end{equation}
converges in the weak-$\ast$ topology to $\Phi$.
\end{proposition}

\begin{proof} Let $\eta = \sum^\oplus\nolimits_{[\zeta]} \zeta$ be the representation mentioned in the statement, in which case $\Phi_n(q) = V^\ast U_n^\ast \eta(q) U_n V$. Then the statement follows from the facts that the limit in~\eqref{Eq:AppEquiv} holds in norm topology of $B(H_\pi)$ and $V$ is a bounded map.\end{proof}

\begin{corollary} The set of u.c.p. maps of the form 
\begin{equation}
\Qq \ni q \mapsto \sum\nolimits_{[\zeta]} V_\zeta^\ast \zeta(q) V_\zeta \in B(H),
\end{equation} 
where the sum runs over a countable subset of $\hat \Qq$ and $V_\zeta$ are bounded linear maps from $H$ to $H_\zeta$, is weak-$\ast$ dense in the space of all u.c.p. maps from $\Qq$ to $B(H)$.
\end{corollary}

\begin{remark}{\rm While weak approximations of states have deficiencies, {\it e.g.} they are not powerful enough to resolve spectral features of the exact state, they serve perfectly well the purpose of computing correlation functions.
}$\Diamond$
\end{remark}

Since this will be our main device for deriving new example of operator product states using Stinespring representations, we present the following example:

\begin{example}{\rm Let $X$ be a compact subset of $\RM$ and $\rho$ a cyclic representation of $C(X)$ on a separable Hilbert space. In the first part of the exercise, we derive for $\rho$ an approximately unitary equivalent representation made out of a countable set of irreducible representations of $C(X)$. In the second part of the exercise, we show how the latter generates a weakly-$\ast$ converging sequence of convex combinations of a fix countable set of pure states of $C(X)$, with the limit state matching to the cyclic representation $\rho$.

 From \cite{DavidsonBook}[Th.~II.1.1], we know that $\rho$ is equivalent to the representation that sends $f \in C(X)$ to the multiplication operator $M_f$ by $f$ on $H_\rho=L^2(X,\mu)$, for some reqular Borel probability measure on $X$. Furthermore, if $T = M_x$ is the multiplication operator by the identity function $f(x) = x$, then $\rho(C(X))$ can be identified with the $C^\ast$-algebra generated by $T$ in $B(H_\rho)$. By Weyl-von Neumann theorem, $T$ is an arbitrarily small compact perturbation of a diagonalizable operator. We need some elements from the proof, which we borrow from \cite{HigsonRoeBook2000}[2.2.5]. Consider a sequence of partitions of $X$ by Borel subsets, such that the diameters of the subsets in the $n$-th partition are equal or less than $\epsilon \, 2^{-n-3}$,  and the $(n+1)$-th partition refines the $n$-th. If $P_n^\epsilon$ denotes the projection onto the finite dimensional linear subspace of $L^2(X,\mu)$ generated by the characteristic functions of the Borel sets of the $n$-th partition,  then $Q_n^\epsilon = P_n^\epsilon -P_{n-1}^\epsilon$ are projections too, $Q_n^\epsilon$ are mutually orthogonal, $\sum_n Q_n^\epsilon$ converges strongly to identity, and $\|Q_n^\epsilon T - T Q_n^\epsilon \| \leq \epsilon 2^{-n }$. Now, 
\begin{equation}\label{Eq:TandQ}
T=\sum_n Q_n^\epsilon T Q_n^\epsilon - \sum_n (Q_n^\epsilon T - T Q_n^\epsilon)Q_n^\epsilon,
\end{equation}
with both sums converging in strong topology, and the first one is a direct sum of  the finite self-adjoint operators $Q_n^\epsilon T Q_n^\epsilon$, while the second one is a compact operator of norm smaller than $\epsilon$. The conclusion is that there exists a diagonal operator $D_\epsilon$ on $\sum^\oplus_n Q_n^\epsilon H_\rho$, a unitary transformation $U_\epsilon : H_\rho \to \sum^\oplus_n Q_n^\epsilon H_\rho$ such that $T = U_\epsilon^\ast D_\epsilon U_\epsilon$ up to a compact operator of norm smaller than or equal to $\epsilon$. Furthermore, it can be arranged to have the spectrum of $D_\epsilon$ contained in the set $X$. Note that the unitary transformation $U_\epsilon$ is a direct sum of finite unitary operators that can be explicitly computed. Now, by Lemma~II.4.3 from \cite{DavidsonBook}, all $D_\epsilon$ operators are approximately unitarily equivalent. Then, by fixing an $\epsilon_0$, there exists a unitary transformation $V_{\epsilon_0,\epsilon} : \sum^\oplus_n Q_n^\epsilon H_\rho \to \sum^\oplus_n Q_n^{\epsilon_0} H_\rho$ such that  $D_\epsilon = V_{\epsilon,\epsilon_0}^\ast D_{\epsilon_0} V_{\epsilon_0,\epsilon}$ plus a compact perturbation whose norm can be taken arbitrarily small, {\it e.g.} smaller than $\epsilon$. By piecing together the parts, we see that  $T = W_{\epsilon_0,\epsilon}^\ast D_{\epsilon_0} W_{\epsilon_0,\epsilon}$ plus a compact correction of norm smaller than or equal to $2 \epsilon$, where $W_{\epsilon_0,\epsilon}= V_{\epsilon_0,\epsilon}U_\epsilon$. Turning now our attention back to the representation $\rho(C(X)) \simeq C^\ast(T)$, we have 
\begin{equation}\label{Eq:DEps}
f(T) - f(W_{\epsilon_0,\epsilon}^\ast D_{\epsilon_0}W_{\epsilon_0,\epsilon})=f(T) - W_{\epsilon_0,\epsilon}^\ast f(D_{\epsilon_0})W_{\epsilon_0,\epsilon} \to 0 \ \ as \ \epsilon \to 0,
\end{equation} 
for any continuous function $f$ over $X$. The irreducible representations of $C(X)$ consist of the evaluations at the points of $X$, $f \mapsto {\rm ev}_x(f) : = f(x)$. If $\{x_k^{\epsilon_0}\}$ are the diagonal entries of $D_{\epsilon_0}$, which is a countable set, then
\begin{equation}
f(D_{\epsilon_0}) \simeq \sum\nolimits_k^\oplus {\rm ev}_{x_k^{\epsilon_0}}(f),
\end{equation}
hence \eqref{Eq:DEps} materializes the statement of Proposition~\ref{Pro:AUEq} in the present context. 

For the second part of the exercise, we proceed as follows:
\begin{equation}
\begin{aligned}
\int d\mu(x) f(x) & = \prescript{}{\mu}{\langle} 1 | M_f | 1 \rangle_{\mu} = \lim_{\epsilon \to 0} \prescript{}{\mu}{\langle} 1 |W_{\epsilon_0,\epsilon} f(D_{\epsilon_0}) W_{\epsilon_0,\epsilon}^\ast |1\rangle_{\mu} \\
& = \lim_{\epsilon \to 0} \sum_k f(x_k^{\epsilon_0}) \prescript{}{\mu}{\langle} 1 |W_{\epsilon_0,\epsilon} P_{\{x_k^{\epsilon_0}\}} W_{\epsilon_0,\epsilon}^\ast|1\rangle_{\mu} .  
\end{aligned}
\end{equation}
where $|1 \rangle_{\mu}$ is the class of the (properly normalized) constant function in $L^2(X,\mu)$ and $P_{\{\cdot\}}$ are the projections corresponding to the specific diagonal entries of $D_{\epsilon_0}$. As promised, the approximate unitary equivalence~\eqref{Eq:DEps} produced a sequence of convex combinations of a fixed set of pure states of $C(X)$, and this sequence converges weakly-$\ast$ to the state corresponding to $\rho$.  
}$\Diamond$
\end{example}

\subsection{Examples produced via Stinespring representations} 
\label{Sec:OpProduct}

The plan is to start from an input data $(\Aa,\Ss,\EM,\phi)$ and apply the reconstruction algorithm, but this time engaging Stinespring representations. To start, we need some preparation. First, we embed $\Ss$ in a $C^\ast$-algebra $\Bb$ and extend the u.c.p. map $\EM$ over $\Aa \otimes \Bb$. Throughout this section, we assume that $\Bb$ is postliminal. Next, we point out that Stinespring's representation engages a positive map with values in $B(H)$ for some Hilbert space $H$. For this, we can compose $\EM$ with a representation of $\Bb$, but, since any representation desintegrates as in Proposition~\ref{Eq:Disinter}, it is fruitful to examine first the families  
\begin{equation}
\EM_\xi : \Aa \otimes \Bb \to B(H_\xi), \quad \EM_\xi : = \xi \circ \EM,
\end{equation}
of u.c.p maps indexed by $[\xi] \in \hat \Bb$. Of course, each $\EM_\xi$ accepts a representation of the form
\begin{equation}\label{Eq:Stinespring1}
\EM_\xi( a \otimes b) = V_\xi^\ast \pi_\xi (a \otimes b) V_\xi \in  B(H_\xi),
\end{equation}
where $V_\xi:H_\xi \to H_{\pi_\xi}$ is an isometry. As the notation suggests, the representation $\pi$ appearing in the above Stinespring representation depends on the chosen representation $\xi$ of $\Bb$. 

\begin{proposition}\label{Prop:RepPi} The representation $\pi_\xi$ is approximately unitarily equivalent to a representation of the type
\begin{equation}\label{Eq:PiRep}
\pi_\xi \sim_a \sum\nolimits^\oplus_{[\zeta]} \pi_{\zeta,\xi} \otimes \zeta,
\end{equation}
where $\pi_{\zeta,\xi}$'s are representations of $\Aa$ and the direct sum seen in the right side is over a countable subset of $\hat \Bb$.
\end{proposition}

\begin{proof} Proposition~\ref{Pro:AUEq} states that every representation of $\Aa \otimes \Bb$ is approximately unitarily equivalent to a direct sum of irreducible representations of $\Aa \otimes \Bb$. Since $\Bb$ is assumed postliminal, then Proposition~\ref{Th:TIrrep} assures us that each irreducible representation of $\Aa \otimes \Bb$ takes the form $\gamma \otimes \zeta$, where both $\gamma$ and $\zeta$ are irreducible representations for $\Aa$ and $\Bb$, respectively. Lastly, any direct sum of terms $\gamma \otimes \zeta$ can be organized over the irreducible representations of $\Bb$, hence it can be brought to the stated form, which incorporates in $\pi_{\zeta,\xi}$ the possible multiplicities over $\zeta$.
\end{proof}

\begin{corollary}\label{Cor:EMxi} Every $\EM_\xi$ can be weakly-$\ast$ approximated by an expression 
\begin{equation}
\tilde \EM_\xi(a \otimes b) = \sum\nolimits_{[\zeta]} V_{\zeta,\xi}^\ast \big(\pi_{\zeta,\xi}(a) \otimes \zeta(b)\big) V_{\zeta,\xi},
\end{equation}
where the sum runs over a countable subset of $\hat \Bb$ and $V_{\zeta,\xi} : H_\xi \to H_{\pi_{\zeta,\xi}} \otimes H_\zeta$ are isometries.
\end{corollary}

Our next task is to resolve the structure of the isometries. For this, it will be convenient to fix basis sets $\{\psi_\mu(\zeta,\xi)\}_\mu$ and $\{\varphi_j(\xi)\}_j$ for the Hilbert spaces $H_{\pi_{\zeta,\xi}}$ and $H_\xi$, respectively.

\begin{proposition}\label{Prop:Iso} Any isometry $V_{\zeta,\xi} :  H_\xi \to H_{\pi_{\zeta,\xi}} \otimes H_\zeta$ can be reduced to the form
\begin{equation}\label{Eq:VForm}
V_{\zeta,\xi} |\varphi(\xi)\rangle  = \sum_{\mu} |\psi_\mu(\zeta,\xi) \rangle \otimes W_\mu(\zeta,\xi)|\varphi(\xi)\rangle,
\end{equation}
where $W_\mu(\zeta,\xi)$'s are bounded linear maps from $H_\xi$ to $H_\zeta$.
Similarly, the adjoint $V_{\zeta,\xi}^\ast :  H_{\pi_{\zeta,\xi}} \otimes H_\zeta \to  H_\xi$ takes the form
\begin{equation}\label{Eq:VCForm}
V_{\zeta,\xi}^\ast |\psi(\zeta,\xi) \rangle \otimes |\varphi(\zeta) \rangle= \sum_{\mu} \langle \psi_\mu(\zeta,\xi)|\psi(\zeta,\xi) \rangle \, W_\mu^\ast(\zeta,\xi)| \varphi(\zeta) \rangle,
\end{equation}
where $W_\mu^\ast(\zeta,\xi)$ is the conjugate of $W_\mu(\zeta,\xi)$. Reciprocally, any family of linear operators $\{W_\mu(\zeta,\xi)\}_\mu$ satisfying the constraints
\begin{equation}\label{Eq:IsoConst}
\sum_\mu W_\mu(\zeta,\xi)^\ast W_\mu(\zeta,\xi) = I_{H_\xi}
\end{equation}
produces an isometry $H_\xi \to H_{\pi_{\zeta,\xi}} \otimes H_\zeta$ via Eq.~\eqref{Eq:VForm}.
\end{proposition}

\begin{proof} Since $V_{\zeta,\xi}$ is a linear bounded map, there exist complex coefficients $\{A_{\mu;i,j}(\zeta,\xi)\}$ such that
\begin{equation}
\begin{aligned}
V_{\zeta,\xi} |\varphi_i(\xi)\rangle & = \sum_{\mu,j} A_{\mu;i,j}(\zeta,\xi)  \, |\psi_\mu(\zeta,\xi) \rangle \otimes |\varphi_j(\zeta)\rangle  \\
& = \sum_{\mu}  |\psi_\mu(\zeta,\xi)\rangle \otimes \sum_j A_{\mu;i,j}(\zeta,\xi) |\varphi_j(\zeta)\rangle.
\end{aligned}
\end{equation}
Then $W_\mu(\zeta,\xi)$ is the linear map from $H_\xi$ to $H_\zeta$ defined by the matrix elements $A_{\mu;i,j}(\zeta,\xi)$ in the obvious basis sets. Furthermore, one can manually check that
\begin{equation}
\big ( |\psi(\zeta,\xi)\rangle \otimes |\varphi(\zeta)\rangle , V_{\zeta,\xi} |\varphi(\xi)\rangle \big )  = \big (  V^\ast_{\zeta,\xi} |\psi_\zeta\rangle \otimes |\varphi_\zeta\rangle , |\varphi_\xi\rangle \big ),
\end{equation} 
if we use the action seen in Eq.~\eqref{Eq:VCForm}. For the last statement, one can check manually that constraint~\eqref{Eq:IsoConst} implies $V_{\zeta,\xi}^\ast V_{\zeta,\xi} = I_{H_\xi}$. Note that this constraint automatically implies that each $W_\mu$ is bounded.
 \end{proof}

\begin{proposition}\label{Pro:TEM} $\EM_\xi$ can be weakly-$\ast$ approximated by an expression of the form
\begin{equation}\label{Eq:EM1}
\tilde \EM_\xi(a \otimes b) = \sum_{[\zeta]} \sum_{\mu,\nu} \langle \psi_\mu(\zeta,\xi)|\pi_{\zeta,\xi}(a)|\psi_{\nu}(\zeta,\xi) \rangle \ W_\mu^\ast(\zeta,\xi)\, \zeta(b) \,  W_\nu (\zeta,\xi). 
\end{equation}
Reciprocally, for any collections of:
\begin{enumerate}[$\ \circ$] 
\item Representations $\{\pi_{\zeta,\xi}, \ [\zeta],[\xi] \in \hat \Bb\}$ of $\Aa$;
\item Basis sets $\{\psi_\mu(\zeta,\xi)\} \subset H_{\pi_{\zeta,\xi}}$;
\item  Bounded maps $\{W_\mu(\zeta,\xi): H_\xi \to H_\zeta, \ [\zeta],[\xi] \in \hat \Bb\}$ satisfying~\eqref{Eq:IsoConst}, 
\end{enumerate}
expression~\eqref{Eq:EM1} defines u.c.p. maps $\tilde \EM_\xi$ from $\Aa \otimes \Bb$ to $B(H_\xi)$.
\end{proposition}

\begin{proof} The direct implication is a simple consequence of Corollary~\ref{Cor:EMxi} and Proposition~\ref{Prop:Iso}. For the reverse implication, if we set $\pi_\xi = \bigoplus_\zeta\,  \pi_{\zeta,\xi} \otimes \zeta$ and $V_\xi = \sum^\oplus_{[\zeta]} V_{\zeta,\xi}$, then expression~\eqref{Eq:EM1} can be assembled back as $V_\xi^\ast \pi_\xi(a \otimes b) V_\xi$, hence a u.c.p. map. 
\end{proof}

\begin{remark}{\rm We want to draw attention to two extremal cases for expression~\eqref{Eq:EM1}. On one end of the complexity spectrum sits the case where only one irreducible representation $\zeta$ is engaged. As we shall see below, iterations of $\tilde \EM$ then involve only this representation $\zeta$. At the other extreme is the case where the sum over $\zeta$ degenerates to an integral $\int_{\hat \Bb} d{\bm \mu}_\xi(\zeta)$ for some regular Borel measure with full support. This case is difficult to iterate and this is the reason we prefer to work with weak-$\ast$ approximants that only involve countable sums.
}$\Diamond$
\end{remark}

\begin{proposition} The system of maps from Proposition~\ref{Pro:TEM} can be iterated:
\begin{equation}\label{Eq:OpProd}
\begin{aligned}
\tilde \EM_{\xi}^{(p)}(\alpha \otimes b)&  = \sum \langle \psi_{\mu_1}(\zeta_1,\xi) |\pi_{\zeta_1,\xi}(a_1) | \psi_{\nu_1}(\zeta_1,\xi) \rangle  \\
& \qquad \cdots 
 \langle \psi_{\mu_p}(\zeta_p,\zeta_{p-1}) |\pi_{\zeta_p,\zeta_{p-1}}(a_p) | \psi_{\nu_p}(\zeta_p,\zeta_{p-1}) \rangle \\
& W_{\mu_1}^\ast(\zeta_1,\xi) \cdots W_{\mu_p}^\ast(\zeta_p,\zeta_{p-1})  \zeta_p(b) W_{\nu_p}(\zeta_p,\zeta_{p-1}) \cdots W_{\nu_1}(\zeta_1,\xi),
\end{aligned}
\end{equation}
where $\alpha=a_1 \otimes \cdots \otimes a_p \in \Aa^{\otimes p}$ and the sum is over all $\zeta$'s, $\mu$'s and $\nu$'s.
\end{proposition}

\begin{proof} The first iteration of $\tilde \EM$ looks as follows: 
\begin{equation}
\begin{aligned}
\tilde \EM_{\xi}\big (a_1 \otimes \tilde \EM(a_2 \otimes b)\big) & = \sum_{[\zeta_1],\mu_1,\nu_1}  \langle \psi_{\mu_1}(\zeta_1,\xi) |\pi_{\zeta_1,\xi}(a_1) | \psi_{\nu_1}(\zeta_1,\xi) \rangle \\ & \qquad \qquad \qquad  W_{\mu_1}^\ast(\zeta_1,\zeta_0)  \zeta_1\big (\tilde \EM(a_2 \otimes b) \big ) W_{\nu_1}(\zeta_1,\xi) .
\end{aligned}
\end{equation}
Since
\begin{equation}
\zeta_1\big (\tilde \EM(a_2 \otimes b) \big )= (\zeta_1 \circ \tilde \EM)(a_2 \otimes b) = \tilde \EM_{\zeta_1}(a_2 \otimes b),
\end{equation}
we can continue
\begin{equation}\label{Eq:StineSbar}
\begin{aligned}
 \tilde \EM_\xi^{(2)}(a_1 \otimes a_2 \otimes b) =&  \sum \langle \psi_{\mu_1}(\zeta_1,\xi) |\pi_{\zeta_1,\zeta_0}(a_1) | \psi_{\nu_1}(\zeta_1,\xi) \rangle \\
 & \qquad \langle \psi_{\mu_2}(\zeta_2,\zeta_1) |\pi_{\zeta_2,\zeta_1}(a_2) | \psi_{\nu_2}(\zeta_2,\zeta_1) \rangle\\
& \qquad \quad W_{\mu_1}^\ast(\zeta_1,\xi) W_{\mu_2}^\ast(\zeta_2,\zeta_1) \zeta_2(b)  W_{\nu_2}(\zeta_2,\zeta_1) W_{\nu_1}(\zeta_1,\xi).
\end{aligned}
\end{equation}
It is then clear that $p$ iterations lead to expression~\eqref{Eq:OpProd}.
\end{proof}

\begin{corollary} The shift map of $\tilde \EM$ on $\Bb$ takes the form
\begin{equation}\label{Eq:SMap}
\begin{aligned}
\xi \circ \bar S^{\circ p}(b)  = \sum  W_{\mu_1}^\ast(\zeta_1,\xi) \cdots & W_{\mu_p}^\ast(\zeta_p,\zeta_{p-1})  \zeta_p(b) \\
& \qquad \qquad W_{\mu_p}(\zeta_p,\zeta_{p-1}) \cdots W_{\mu_1}(\zeta_1,\xi).
\end{aligned}
\end{equation}
\end{corollary}

\begin{remark}{\rm The conditions of Th.~\ref{Th:Main3} are satisfied if $\xi \circ \bar S^{\circ p}$ converges to $1_{B(H_\xi)}$ as $p \to \infty$, for each $[\xi] \in \hat \Bb$. The right side of~\eqref{Eq:StineSbar} should help us decide if that is the case or not.
}$\Diamond$
\end{remark}

Our last task is to evaluate the constructed state: 
\begin{equation}
\omega(\alpha) = \big ( \phi \circ \EM^{(p)}\big )(\alpha \otimes 1), \quad \alpha \in \Aa^{\otimes p}, \quad p \geq 1.
\end{equation}
By embedding $\Bb$ into $\sum_\xi^\otimes B(H_\xi)$, we finally have:

\begin{proposition} Every shift-invariant state over $\Aa_\ZM$ can be weakly-$\ast$ approximated as
\begin{equation}
\Aa^{\otimes p} \ni \alpha \mapsto \tilde \omega(\alpha) = \sum_\xi \big (\phi_\xi \circ \tilde \EM_\xi^{(p)}\big )(\alpha \otimes 1),
\end{equation}
where $\{\phi_\xi\}$ is a system of c.p. maps on $\{B(H_\xi)\}$, $\xi \in \hat \Bb$.
\end{proposition}

\begin{example}\label{Ex:HilbertC}{\rm Consider $\Aa = C([0,1])$, the algebra of continuous functions over the closed interval $[0,1]$. Then $\Aa_\ZM$ coincides with $C\big ([0,1]^{\times \ZM}\big )$, the $C^\ast$-algebra of continuous functions over the Hilbert cube. The shift-invariant states over $\Aa_\ZM$ are one on one with the shift-invariant normalized Radon measures over $[0,1]^{\times \ZM}$. Our theory says that ergodic measures can be generated via the reconstruction process from a reduce data $\big ( \Aa = C([0,1]),\Ss,\EM,\phi)$, provided we can verify the conditions of Theorem~\ref{Th:Main3}. We supply here an example of such data. For this, we take the embedding $C^\ast$-algebra of $\Ss$ to be $M_2(\CM)$, whose spectrum consists of a single point, corresponding to the identity representation. Thus, the sum over $\zeta$ in expression~\eqref{Eq:EM1} disappears and we only need to fix a single representation $\pi$ of $\Aa$. If $\PM$ is a probability measure over the interval $[0,1]$, we take $\pi$ to be the GNS representation corresponding to the state $a \mapsto \int {\rm d}\PM(t) \, a(t)$. The Hilbert space of this representation is $H_\pi=L^2([0,1],{\rm d}\PM)$ and the functions from $\Aa$ act by multiplication on the square integrable functions from $L^2([0,1],{\rm d}\PM)$. We fix a basis $\{\psi_\mu\}$ for $H_\pi$ and, for $\mu =0,1,2,3$ and $d < 1/\sqrt{3}$, we pick
\begin{equation}
W_0= \sqrt{1-3d^2}\, \sigma_0, \quad W_\mu = d \, \sigma_\mu, \quad \mu=1,2,3,
\end{equation}
where $\sigma_0$ is the unit matrix and the rest of the $\sigma$'s are Pauli's matrices. We set $W_\mu =0$ for the remaining values of $\mu$. With these choices, $\sum_{\mu} W_\mu^\ast W_\mu = \sigma_0$ and, as such, the map
\begin{equation}\label{Eq:EMx}
\EM(a \otimes b) = \sum_{\mu,\nu=0}^3 \int {\rm d}\PM(t) \, \bar \psi_\mu(t)\psi_{\nu}(t) a(t) \ W_\mu^\ast\, b \,  W_\nu
\end{equation}
is u.c.p.. Furthermore, if we express $b \in M_2(\CM)$ as $b = b_0 \sigma_0 + \vec b \cdot \vec \sigma$, we have
\begin{equation}\label{Eq:EMy}
\bar S(b) = \EM(1 \otimes b) = \sum_{\mu=0}^3  W_\mu^\ast\, b \,  W_\mu = b_0 \sigma_0 + \lambda \vec b \cdot \vec \sigma,
\end{equation}
where $\lambda = (1-4d^2)$ is a number of absolute value strictly smaller than 1. Therefore, $\bar S^{\circ n}(b)$ converges as $n \to \infty$ to ${\rm tr}(b) \sigma_0$, where ${\rm tr}$ is the unique trace state over $M_2(\CM)$. As a consequence, we are in the conditions of Th.~\ref{Th:Main3}, hence we know for sure that the reduced space of the constructed state in $M_2(\CM)$ and, furthermore, that the only choice for a $\bar S$-invariant state is $\phi(b) = {\rm tr}(b)$. Lastly, if $\beta$ is the embedding of $a_1 \otimes 1^{\otimes p} \otimes a_2\in \Aa^{\otimes (p+2)}$ in $\Aa_R$, we can compute the correlation function
\begin{equation}
\omega(\beta) = ({\rm tr} \circ \EM)(a_1 \otimes \bar S^{\circ p}(\EM(a_2\otimes e)))
\end{equation}
explicitly as
\begin{equation}
\begin{aligned}
\omega(\beta) =  \omega(\beta_1)\omega(\beta_2) + \lambda^p \sum_{\mu,\nu} & \int {\rm d}\PM(t)  \, \bar \psi_{\mu_1}(t)\psi_{\nu_1}(t) a_1(t) \\
& \quad \int {\rm d}\PM(t)  \, \bar \psi_{\mu_2}(t)\psi_{\nu_2}(t) (a_2(t)-\omega(a_2)) \\
&  \qquad \quad {\rm tr}(W_{\mu_1}^\ast\, W_{\mu_2}^\ast W_{\nu_2} \,  W_{\nu_1}),
\end{aligned}
\end{equation}
where $\beta_{1,2}$ are (any) embeddings of $a_{1,2}$ in $\Aa_\ZM$. To conclude, we have constructed a Radon measure $\Gamma$ over the Hilbert cube that has a known correlation decay law, in the sense that
\begin{equation}
\lambda^{-p}\Big [\int {\rm d}\Gamma(\{t_j\}) a_1(t_n)a_2(t_{n+p}) - \int {\rm d}\Gamma(\{t_j\}) a_1(t_n)\int {\rm d}\Gamma(\{t_j\}) a_2(t_{n+p}) \Big ]
\end{equation} 
converges to a constant for $p \to \infty$. 

This example demonstrates that, even though we are dealing with a commutative $C^\ast$-algebra $\Aa_\ZM$, the reduced space $\Bb_\omega$ does not necessarily have to be embedded into a commutative $C^\ast$-algebra.\footnote{Recall that $\Bb_\omega$ has only a linear structure, hence no statement can be made about a multiplicative structures before embedding.} If we replace the interval $[0,1]$ by a finite set of points, it is known that $\Ss$ can be always taken a finite commutative $C^\ast$-algebra (see example 7.1 in \cite{FannesCMP1992}). Things are different for the present context because, in general, we cannot choose a basis for $L^2([0,1],{\rm d}\PM)$ such that all matrix elements seen in \eqref{Eq:OpProd} are diagonal. 
}$\Diamond$
\end{example}

\end{document}